\documentclass[11pt]{article}
\usepackage{amsmath}
\normalbaselines

\renewcommand{\Box}{\mbox{}}

\newtheorem{lem}{Lemma}[section]
\newtheorem{prop}{Proposition}[section]

\begin{document}
\bibliographystyle{unsrt}

\def\bea*{\begin{eqnarray*}}
\def\eea*{\end{eqnarray*}}
\def\ba{\begin{array}}
\def\ea{\end{array}}
\count1=1
\def\be{\ifnum \count1=0 $$ \else \begin{equation}\fi}
\def\ee{\ifnum\count1=0 $$ \else \end{equation}\fi}
\def\ele(#1){\ifnum\count1=0 \eqno({\bf #1}) $$ \else \label{#1}\end{equation}\fi}
\def\req(#1){\ifnum\count1=0 {\bf #1}\else \ref{#1}\fi}
\def\bea(#1){\ifnum \count1=0   $$ \begin{array}{#1}
\else \begin{equation} \begin{array}{#1} \fi}
\def\eea{\ifnum \count1=0 \end{array} $$
\else  \end{array}\end{equation}\fi}
\def\elea(#1){\ifnum \count1=0 \end{array}\label{#1}\eqno({\bf #1}) $$
\else\end{array}\label{#1}\end{equation}\fi}
\def\cit(#1){
\ifnum\count1=0 {\bf #1} \cite{#1} \else 
\cite{#1}\fi}
\def\bibit(#1){\ifnum\count1=0 \bibitem{#1} [#1    ] \else \bibitem{#1}\fi}
\def\ds{\displaystyle}
\def\hb{\hfill\break}
\def\comment#1{\hb {***** {\em #1} *****}\hb }

\newcommand{\TZ}{\hbox{T\hspace{-5pt}T}}
\newcommand{\MZ}{\hbox{I\hspace{-2pt}M}}
\newcommand{\ZZ}{\hbox{Z\hspace{-3pt}Z}}
\newcommand{\NZ}{\hbox{I\hspace{-2pt}N}}
\newcommand{\RZ}{\hbox{I\hspace{-2pt}R}}
\newcommand{\CZ}{\,\hbox{I\hspace{-6pt}C}}
\newcommand{\PZ}{\hbox{I\hspace{-2pt}P}}
\newcommand{\QZ}{\hbox{I\hspace{-6pt}Q}}
\newcommand{\KZ}{\hbox{I\hspace{-2pt}K}}

\vbox{\vspace{38mm}}
\begin{center}
{\LARGE \bf  Toric Representation and Positive Cone of Picard Group and  Deformation Space in Mirror Symmetry of Calabi-Yau  Hypersurfaces in Toric Varieties }\\[5mm]

Shi-shyr Roan
\\{\it Institute of Mathematics \\ Academia Sinica \\ 
Taipei , Taiwan \\ (e-mail: maroan@gate.sinica.edu.tw)} \\[5mm]
\end{center}

\begin{abstract} We derive the combinatorial representations of Picard group and  deformation space of anti-canonical hypersurfaces of a toric variety using techniques in toric geometry. The mirror cohomology correspondence in the context of mirror symmetry is established for a pair of Calabi-Yau (CY) ${\sf n}$-spaces in toric varieties defined by reflexive polytopes for an arbitrary dimension ${\sf n}$. We further identify the Kahler cone of the toric variety and degeneration cone of CY hypersurfaces, by which the Kahler cone and degeneration cone for a mirror CY pair are interchangeable under mirror symmetry. In particular, different degeneration cones of a CY 3-fold are corresponding to flops of its mirror 3-fold.

\end{abstract}

\par \vspace{5mm} \noindent
 \par \noindent
{\rm 1999 PACS}: 02.10.Rn, 02.10.Eb   \par \noindent

\section{Introduction}
\setcounter{equation}{0}
For a quasi-smooth Calabi-Yau (CY) ${\sf n}$-space $X$ in a toric variety, the Picard group ${\rm Pic}(X)_{\CZ}$ and deformation space ${\rm Def}(X)_{\CZ}$ in the Hodge spaces, $H^{1,1}(X)$ and $H^{{\sf n}-1,1}(X)$ respectively, play an important role in the algebraic geometry study of $X$ and  its applications in string physics (see, e.g. \cite{GY, Y} and references therein). In the case when the toric variety is defined by a simplicial-cone decomposition with certain triangulation conditions, the anti-canonical hypersurface $X$ is indeed a smooth manifold for ${\sf n} \leq 3$. In particular for the case ${\sf n}=2$,   
\bea(ll)
X= {\rm K3 ~ surface}: & H^{1,1}(X) = {\rm Pic}(X)_{\CZ} \oplus {\rm Def}(X)_{\CZ}.
\elea(K3)
For ${\sf n}$ greater than 2, the Picard group and deformation space are the same as the Hodge spaces: 
\bea(lll)
 H^{1,1}(X) = {\rm Pic}(X)_{\CZ}, & H^{{\sf n}-1,1}(X)  = {\rm Def}(X)_{\CZ} & {\rm for} ~ {\sf n} \geq 3 .
\elea(PiDe)
By the mirror symmetry of Calabi-Yau (CY) spaces, we mean a pair of quasi-smooth CY ${\sf n}$-spaces, $X$ and $X^*$, with a canonical identification of Picard group and deformation space:
\bea(ll)
{\rm Pic}(X)_{\CZ} \simeq {\rm Def}(X^*)_{\CZ} , &  {\rm Def}(X)_{\CZ} \simeq  {\rm Pic}(X^*)_{\CZ},
\elea(PiDe)
which are compatible with (quantum) cohomology ${\sf n}$-product. For ${\sf n}=2$, the mirror symmetry (\req(PiDe)) signifies the interchangeable relation between the algebraic cycles and the transcendental cycles of K3 surfaces $X$ and $X^*$. This symmetry is linked to the Arnold duality for the 14 exceptional singularities of modality one \cite{Ko, R92}. For ${\sf n}\geq 3$, by (\req(PiDe)), the mirror symmetry (\req(PiDe)) implies 
\be
h^{1, 1}(X) = h^{{\sf n}-1,1}(X^*) ~ ~ ~ ~ ({\sf n}\geq 3 ).
\ele(h=h)
In particular when ${\sf n}=3$, the Euler numbers of $X$ and $X^*$ are of opposite sign, $\chi (X) = - \chi (X^*)$. 
Indeed, the mirror symmetry was first found in \cite{GP, Y} as a pair of Fermat-type CY 3-folds corresponding to the same conformal field theory with a reversal of left $U(1)$-charge in string theory.  Subsequently, the mirror cohomology correspondence (\req(PiDe)) for Fermat-type CY 3-folds was shown in \cite{R91} by the method of toric geometry as anti-canonical hypersurfaces of abelian quotients of weighted projective 4-spaces. A similar discussion for higher dimensional Fermat-type CY spaces was also given in \cite{R93y}. The mirror symmetry was further extended to a large class of CY spaces in \cite{B}, where the anti-canonical hypersurfaces in $n$-dimensional toric varieties defined by reflexive polytope and its dual polytope for $n (= {\sf n}+1)$ are found to satisfy the Hodge-number equality (\req(h=h)) (\cite{B} Theorem 4.4.3). However, to the best of the author's knowledge, the mirror correspondence (\req(PiDe)) has not been established in literature till now,  except the simplex polytope (i.e. Fermat-type) case for dimension $n=4$ \cite{R91}. The importance of a canonical identification of toric representatives in (\req(PiDe)) is stressed 
as the keys either to  string conformal-field applications, or to enumerative geometry about rational-curve counting in a CY space.  With a correct cohomology identification in (\req(PiDe)), it is believed that the quantum cohomology computation of a CY manifold can be carried out by the method of complex-structure variation in its mirror degeneration family, as suggested by works in \cite{COGP, HKLY, Y} about the rational-curve problem in Fermat quintic or some other specific CY 3-folds. The object of this paper is to establish the canonical isomorphism (\req(PiDe)) between  a mirror pair of quasi-smooth CY hypersurfaces in $n$-dimensional toric varieties defined by reflexive polytopes for  $n \geq 2$. In this work, we employ techniques in toric geometry to build a special kind of combinatorial representations for Picard group and deformation space of those CY spaces, by which the mirror correspondence is established by identifying toric representatives of the representation.  Of particular interest is the effect of the topology of CY spaces on those representations in mirror symmetry. It is known that the cohomology product of Picard group in (\req(PiDe)) depends on the topology of manifolds, even in the Fermat CY case. Indeed, the Fermat CY 3-folds in \cite{GP, R91} are constructed as crepant resolutions of a (singular) hypersurface of an abelian-quotient of weighted projective $4$-space \cite{MOP, R89, RY}.  Such crepant  resolutions are not unique in general. Two different ones are connected by a process of flops in birational geometry of 3-folds  with different topological triple-products of Picard group (see, e.g. \cite{R93}), upon which the quantum product builds, albeit they both have the same cohomology representation. In the mirror correspondence (\req(PiDe)), the effect of flop of a CY 3-fold on the deformation of its mirror  has not been previously investigated in the literature, even for the Fermat-type case \cite{R00}. Indeed, it is one main purpose of this present work to search an equivalent notion on the deformation space in mirror symmetry which corresponds to the flop in birational geometry. After sorting out the toric data between mirror Fermat-type CY 3-folds, we find that the topological change due to flop of a CY $3$-fold is linked to a different path of maximal unipotent degeneration in the moduli space of its mirror CY spaces. Consequently,  by using techniques in toric geometry, one expects the same conclusion should also be valid for a mirror pair of anti-canonical hypersurfaces in toric varieties defined by reflexive polytope and its dual polytope. In the present paper we show that it is indeed the case. However,  it turns out that the identifications in (\req(PiDe)) are quite involved, especially for the positive cone structure. Various technical details need to be checked before we then explicitly verify our proposal in a precise form.

This paper is organized as follows. Section \ref{sec.PG} is devoted to structures related to Picard group of hypersurfaces of a toric variety.  We start with the homogeneous coordinate system of a toric variety in Subsection \ref{ssec.Cood}. In Subsection \ref{ssec.Pt}, we recall some facts in \cite{R96} about the Picard group of a generic anti-canonical hypersurface of a toric variety defined by reflexive polytope that will be used in the following. In Subsection \ref{ssec.Kcone}, we present an effective mechanism to extract Kahler cone of a toric variety from its rational simplicial-cone toric structure. Explicit calculation of some examples are performed here.
In Section \ref{sec.Defm}, we investigate structures related to the deformation space of anti-canonical hypersurfaces of a toric variety. In the discussion of this section, we use some facts in toric geometry  of which some detailed argument with technical complexity are provided in Appendix for easy reference.  We first recall the result in \cite{R96} about a combinatorial basis of anti-canonical sections of a toric variety, then describe the coordinate form of the basis elements. In Subsection \ref{ssec.DX0}, we determine the deformation classes of anti-canonical hypersurfaces in the "minimal" toric variety. First, we find the moduli space of the anti-canonical hypersurface by Jacobian-ring technique, then identify  the deformation equivalent classes modulus the transformation of toric-variety automophisms.  By the cohomology relation of hypersurfaces between a toric variety and the minimal toric variety, we derive the combinatorial representation of the deformation space of an anti-canonical hypersurface in toric variety defined by reflexive polytope in Subsection \ref{ssec.DX}. In Subsection \ref{ssec.+deg}, we define the degeneration cone in the moduli of CY hypersurfaces, a structure parallel to the Kahler cone of the toric variety in Subsection \ref{ssec.Kcone}. Using the combinatorial representation of Picard group and deformation space established in Sections \ref{sec.PG} and \ref{sec.Defm}, we show the mirror correspondence (\req(PiDe)) and the identification of Kahler and degeneration cones between the mirror anti-canonical hypersurfaces of toric varieties defined by a reflexive and its dual polytopes in Section \ref{sec.mir}.
We close in Section \ref{sec.F} with concluding remarks. Some results in \cite{R96} are reviewed in Appendix, where some basic facts with technical arguments in toric geometry are also presented  for easy use and reference of this paper (for the details, also see e.g. \cite{Dan, KKMS, Oda}).
\par \vspace{0.2in} \noindent
{\bf Notation.} 
In this work, we use the following notations. Let $L$ be a $n$-dimensional lattice $(\simeq \ZZ^n)$, and $L^* (= \mbox{Hom} ( L , \ZZ ) ) $ the dual lattice of $L$.  We set $L_{\KZ} = L \bigotimes_{\ZZ} {\KZ} $  for ${\KZ} = {\QZ} , {\RZ}$ or ${\CZ}$. Denote ${\bf T} ( L ) =  L_{\CZ}  / L  ~ ( \simeq {\CZ}^{* n} )$  the (algebraic) $n$-torus whose 1-parameter subgroups and characters are identified with elements in $L$ and $L^*$  respectively. For a complete rational fan $\Sigma$ in $L_{\RZ}$, i.e. a polyhedral cone decomposition of $L_{\QZ} \subseteq L_{\RZ}$, 
we shall denote ${\PZ}_{\Sigma} (={\PZ}_{(\Sigma, L)})$ the complete ${\bf T} ( L )$-toric variety  associated to $\Sigma$ \cite{Dan, KKMS, Oda}. In general, ${\PZ}_{\Sigma}$ is a singular space, with at most abelian quotient singularities when  $\Sigma$  is a simplicial cone decomposition. By an integral polytope $\triangle \subseteq L_{\RZ}$, we mean a $n$-dimensional convex hull generated by finitely many lattice points of $L$, whose interior contains the origin of $L$, and the set of vertices of $\triangle$ will be denoted by ${\cal V} ( \triangle )$. A $(n-1)$-dimensional  face will also be called a facet of $\triangle$. The face-decomposition of the boundary $\partial \triangle$ of a polytope $\triangle$  gives rise to the polyhedral-cone decomposition $\Sigma_0$ in $L_{\RZ}$. The ${\bf T} ( L )$-toric variety ${\PZ}_{\Sigma_0}$ defined by $\Sigma_0$ is minimal in the sense that  ${\PZ}_{\Sigma_0}$ is dominated by every ${\bf T} ( L )$-toric variety ${\PZ}_{\Sigma}$ for a refinement $\Sigma$ of $\Sigma_0$, (equivalently, $\Sigma \cap \partial \triangle$ is polytope decomposition of $\partial \triangle$). In particular, a simplicial decomposition  $\Lambda$ of $\partial \triangle$ with $\Lambda^{(0)} \subseteq \partial \triangle \cap L$, where $\Lambda^{(j)}$ denotes the collection of $j$-simplices of $\Lambda$ for $0 \leq j \leq n-1$, gives rise to a complete fan $\Sigma ( = \Sigma (\Lambda))$ of $L_{\RZ}$ consisting of cones $\sigma ( {\tt s} ) (= \sum_{\delta \in {\tt s} \cap \Lambda^{(0)}} \RZ_{\geq 0} \delta )$ for simplices $ {\tt s} \in \Lambda $. Then  ${\PZ}_{\Sigma} (= {\PZ}_{\Sigma (\Lambda)})$ is a complete ${\bf T} ( L )$-toric variety dominating ${\PZ}_{\Sigma_0}$  with at most abelian quotient singularities.

Denote $\langle * , \ * \rangle :  L^*_{\RZ} \times L_{\RZ} \longrightarrow \RZ $ the canonical pairing with integral values on $L^* \times L$. 
The dual polytope $\triangle^*$ of $\triangle$ in $L^*_{\RZ}$ is defined by
$$ 
\triangle^*  = \{ y \in  L^*_{\RZ} \ | \ \langle  y , x \rangle \ \geq -1 \ 
, \ \mbox{for } \ x \in \triangle \ (\Leftrightarrow x \in {\cal V} ( \triangle ) ) \}. 
$$
The dual face $F^*$ of a $m$-dimensional polyhedral face $F$ of $\triangle  ~ (0 \leq m \leq n-1 )$ is  the $(n-m-1)$-dimensional face of $\triangle^*$ defined by
$$
F^* =  \{ y \in \triangle^* \ | \ \langle  y, x \rangle = -1 
\mbox{ \ for \ } x \in F \ (\Leftrightarrow x \in F \cap {\cal V} ( \triangle )  \} .
$$
\par \noindent
{\bf Definition.} $( \triangle , L )$ is a reflexive polytope  iff 
both $(\triangle, L)$ and $(\triangle^*, L^*)$ are integral \cite{B}. \par \vspace{0.1in}  \noindent
It is known that the origin is the only lattice point in the interior $\mbox{ Int } ( \triangle )$ of a reflexive polytope $( \triangle , L )$. For a reflexive polytope $(\triangle, L)$, a simplicial decomposition  $\Lambda$ of $\partial \triangle$ will always be assumed to satisfy the following condition
\be
\Lambda^{(0)} = L \ \bigcap \ \big( \partial \triangle - 
\bigcup \{ \mbox{Int}(F)  \ | \  F :  \mbox{codim-1 \ face \
of \ } \triangle \} \big). 
\ele(Lam0)
Note that the toric variety ${\PZ}_{\Sigma}$ is Gorenstein over ${\PZ}_{\Sigma_0}$. Throughout this paper except in Appendix or otherwise stated, a polytope $(\triangle, L )$ will always be a reflexive polytope with a simplicial decomposition $\Lambda$ of $\partial \triangle$ satisfying (\req(Lam0))
and the relation\footnote{For simple notations, we assume the relation (\req(Lam0)) here, without which the conclusions of this paper are still valid as argued in \cite{R96} by regarding ${\PZ}_{\Sigma (\Lambda)}$ as a $(L/L_0)$-quotient of ${\PZ}_{\Sigma (\Lambda_0)}$ with the homogeneous coordinates in Section \ref{ssec.Cood}, and identifying the anti-canonical hypersurfaces of ${\PZ}_{\Sigma (\Lambda)}$ with $(L/L_0)$-invariant anti-canonical hypersurfaces of ${\PZ}_{\Sigma (\Lambda_0)}$.}
\be
L =  {\rm the ~ sublattice ~} L_0 ~ {\rm of ~ } L ~ {\rm generated ~ by ~ } \Lambda^{(0)}. 
\ele(L0) 
The irreducible toric divisors in ${\PZ}_{\Sigma (\Lambda)}$ are determined by 1-skeleton of $\Sigma (\Lambda)$, parametrized by $\Lambda^{(0)}$. Denote by $e^\delta$ the toric divisor in ${\PZ}_{\Sigma (\Lambda)}$ associated to an element $\delta \in \Lambda^{(0)}$. The divisor lattice is defined by 
\bea(lll)
D_\triangle (= D_{(\triangle, L)}) := \bigoplus_{\delta \in \Lambda^{(0)}} \ZZ e^\delta & ( \simeq \ZZ^d ) , & d :=|\Lambda^{(0)}|.
\elea(DL)
By (\req(L0)), we have the following exact sequence of abelian groups:
\bea(ll)
0 \longrightarrow {\bf n}_\triangle \stackrel{\iota}{\longrightarrow} D_\triangle \stackrel{\beta}{\longrightarrow} L \longrightarrow 0, & \beta ( e^\delta ) := \delta ,
\elea(exa)
where ${\bf n}_\triangle   ( = {\bf n}_{(\triangle, L)}) : =  \mbox{Ker} (\beta)$ is the rank-$(d-n)$ sublattice of $D_\triangle $  with $\iota$ the inclusion morphism. The dual sequence of (\req(exa)) is
\bea(ll)
0 \longrightarrow L^* \stackrel{\beta^*}{\longrightarrow} D_\triangle^\dagger
\stackrel{\iota^*}{\longrightarrow} {\bf n}_\triangle^\dagger \longrightarrow  0  , & \bigg( D_\triangle^\dagger :={\rm Hom}(D_\triangle, \ZZ) , ~ {\bf n}_\triangle^\dagger :={\rm Hom}({\bf n}_\triangle , \ZZ) \bigg),
\elea(exa*)
where the basis of $D_\triangle^\dagger$ dual to $e_{\delta}$'s are denoted by 
\be
D_\triangle^\dagger = \bigoplus_{\delta \in \Lambda^{(0)}} \ZZ e^{\delta \dagger} , ~ ~ \langle e^{\delta^\dagger}, e^{\delta'} \rangle = \varepsilon_{\delta}^{\delta'} ~ ( :=1 ~ {\rm if} ~ \delta = \delta', ~ ~ 0 ~ {\rm otherwise}).
\ele(D*)
For a cone $C$ in a $\RZ$-vector space $V$, the dual cone of $C$ in  the dual space will be denoted by 
$$
\widehat{C} = \{ x \in {\rm Hom}(V, \RZ) | \ \langle x , v \rangle \geq 0 ~ ( v \in C) \}. 
$$

\section{Picard Group of Anti-canonical Hypersurface of Toric Variety\label{sec.PG}}
\setcounter{equation}{0}
We start with coordinate systems of  toric varieties for later use, then review  
some facts in \cite{R96} about the Picard group of a generic anti-canonical hypersurface of a toric variety defined by reflexive polytope. This will also serve to establish the notation. The Kahler cone of toric variety will be discussed in the last subsection here. 

\subsection{Homogeneous coordinates of a toric variety \label{ssec.Cood}} 
First, we describe the homogeneous coordinate system of the toric variety ${\PZ}_{\Sigma (\Lambda)}$ in \cite{Au, Cox, O94, R96}. Regard (\req(exa)) and (\req(exa*)) as the 1-parameter subgroups and characters of the following exact sequence of tori:
$$
0 \longrightarrow {\bf T} ( {\bf n}_\triangle ) \longrightarrow {\bf T} ( D_\triangle) 
\longrightarrow {\bf T} ( L ) \longrightarrow 0 .
$$
Denote
$$
D_{\triangle~\KZ} = \bigoplus_{\delta \in \Lambda^{(0)}} \KZ e^\delta ~ ( \simeq  \KZ^d ) ~ ~ {\rm for} ~ \KZ = \QZ, \RZ, \CZ ,
$$
and the first quadrant of $D_{\triangle~\RZ }$ by
$$
\Omega (= \Omega_\triangle) := \sum_{\delta \in \Lambda^{(0)}} \RZ_{\geq 0} e^\delta .
$$
The ${\bf T} ( D_\triangle)$-toric variety associated to the face-decomposition of $\Omega$ is $D_{\triangle~\CZ} $ with coordinates $\sum_{\delta} z_\delta e^\delta$ for $z_\delta$ corresponding to $e^{\delta \dagger}$ in (\req(D*)):
$$
{\bf T} ( D_\triangle) \simeq \prod_{\delta \in \Lambda^{(0)}} \CZ^* \subseteq D_{\triangle~\CZ}  \ni z = (z_\delta)_{\delta \in \Lambda^{(0)}} .
$$ 
Associated to a triangulation $\Lambda$ of $\partial \triangle$, there is the (integral) simplicial fan $\widetilde{ \Sigma} \subseteq \partial \Omega$ in $D_{\triangle~\RZ}$ lying over the fan $\Sigma (\Lambda)$ of $L_{\RZ}$,
$$
\widetilde{\Sigma} ~ (= \widetilde{\Sigma} (\Lambda)) = \{ \widetilde{\sigma}({\tt s}) | {\tt s}: {\rm simplex ~ in } ~ \Lambda \}, ~ ~ ~ \widetilde{\sigma}({\tt s}): = \sum _{\delta \in {\tt s} \cap \Lambda^{(0)}} \RZ_{\geq 0} e^\delta \subset D_{\triangle~\RZ}.
$$
Then $\widetilde{\Sigma}$ gives rise to a ${\bf T} ( D_\triangle)$-toric variety, denoted by $\CZ_{\widetilde{\Sigma} } (= \CZ_{\widetilde{\Sigma} (\Lambda)} )$, which is an affine open subset of $D_{\triangle~\CZ}$  given by
\be
{\bf T} ( D_\triangle) \subset \CZ_{\widetilde{\Sigma} } = D_{\triangle~\CZ} -  \bigcup_{I } \{ \sum z_{\delta} e^\delta
\ | \ z_{\delta} = 0 \mbox{ for } \ \delta \in I  \} 
\ele(CZ)
where the index $I$ runs over subsets of $\Lambda^{(0)}$ which are not in the form ${\tt s} \cap \Lambda^{(0)}$ for some 
${\tt s} \in \Lambda$. Since $\beta$ in (\req(exa)) induces a fan-correspondence from 
$\widetilde{\Sigma }$ to $\Sigma$, sending $\widetilde{\sigma}$ to $\sigma $, it gives rise to the "principal ${\bf T} ( {\bf n}_\triangle )$-bundle"  
\be
\pi : \CZ_{\widetilde{\Sigma}} \longrightarrow {\PZ}_{\Sigma}  ~ ~ ( = \CZ_{\widetilde{\Sigma}}/{\bf T} ( {\bf n}_\triangle )) .
\ele(prin)
The coordinates $z = (z_\delta)_{\delta }$ of $\CZ_{\widetilde{\Sigma}}$ will be called the homogeneous coordinates of ${\PZ}_{\Sigma}$. Similarly, there is the coordinate system of the minimal toric variety ${\PZ}_{\Sigma_0}$ defined by the face-decomposition of $\partial \triangle$, by considering the following exact sequences as (\req(exa)) and (\req(exa*)):
\bea(lll)
0 \longrightarrow {\bf n}_{\triangle ~ 0} \stackrel{\iota_0}{\longrightarrow} D_{\triangle ~ 0} \stackrel{\beta_0}{\longrightarrow} L \longrightarrow 0, & ~ &
0 \longrightarrow  L^* \stackrel{ \beta_0^*}{\longrightarrow } D_{\triangle ~ 0}^\dagger
\stackrel{\iota_0^* }{\longrightarrow }  {\bf n}_{\triangle ~ 0}^\dagger \longrightarrow   0  
\elea(exa0)
where $D_{\triangle~0} := \bigoplus_{\upsilon \in {\cal V} ( \triangle )} \ZZ e^\upsilon$ , $\beta_0 ( e^\upsilon ) := \upsilon$, and ${\bf n}_{\triangle 0}  : =  \mbox{Ker} (\beta_0)$ .  As in (\req(prin)), there is the homogeneous coordinates, $\zeta = (\zeta_\upsilon)_{\upsilon \in {\cal V} ( \triangle )}$, of ${\PZ}_{\Sigma_0}$ with the "${\bf T} ( {\bf n}_{\triangle 0})$-principal bundle" 
\bea(l)
\pi_0 : \CZ_{\widetilde{\Sigma}_0} \longrightarrow {\PZ}_{\Sigma_0},
\elea(prin0)
where $\CZ_{\widetilde{\Sigma}_0}$  is the complement of $\bigcup_{I_0 \notin \{ F \cap {\cal V} ( \triangle ) | F : {\rm face ~ of} ~ \triangle \} } \{ \sum \zeta_\upsilon e^\upsilon \ | \ \zeta_\upsilon = 0 \mbox{ for } \ \upsilon \in I_0  \}$ in $D_{\triangle ~ 0 ~\CZ}$. The fibration (\req(prin)) and (\req(prin0)) are related by the following morphisms between $D_{\triangle ~ 0}$ and $D_{\triangle}$:
\bea(lll)
j: D_{\triangle ~ 0} \hookrightarrow D_{\triangle}, & j (e^\upsilon) = e^\upsilon  &{\rm for} ~ \upsilon \in {\cal V} ( \triangle ) , \\
 p : D_\triangle \rightarrow D_{\triangle ~ 0}, & p (e^\delta) = \sum_{\upsilon \in F \cap {\cal V} ( \triangle )} r^{F,\delta}_{ \upsilon} e^\upsilon  & {\rm for} ~ \delta \in   \Lambda^{(0)} ,
\elea(pDD)
where $F$ is a face of $\triangle$ whose interior ${\rm Int}(F) $ contains $\delta$, and $r^{F,\delta}_{\upsilon}$'s are the $F$-face expression of $\delta$, i.e.  $\delta = \sum_{\upsilon \in F \cap {\cal V} ( \triangle )} r^{F,\delta}_{\upsilon} \upsilon$ for  $r^{F,\delta}_{\upsilon} \in \QZ_{>0}$  with $\sum_{\upsilon' \in F \cap {\cal V} ( \triangle )} r^{F,\delta}_{\upsilon} =1$. Then 
$$
p \cdot j = {\rm id}_{D_{\triangle~0}}, ~ ~ ~ \beta \cdot j = \beta_0 , ~ ~ ~ \beta_0 \cdot p = \beta ,
$$
where $\beta, \beta_0$ are morphisms in (\req(exa)) and (\req(exa0)) respectively.  Then $j, p$ induce an embedding and projection between ${\bf n}_{\triangle ~ 0}$ and ${\bf n}_{\triangle}$. The morphism
\be
p^*: D_{\triangle ~ 0}^\dagger \longrightarrow D_\triangle^\dagger,~ ~ p^*(e^{\upsilon \dagger})= \sum_{\upsilon \in F, \delta \in {\rm Int}(F)  \cap \Lambda^{(0)}} r^{F, \delta}_{\upsilon } e^{\delta \dagger} ~ {\rm for} ~ v \in {\cal V} ( \triangle )
\ele(p*)
induced by $p$ in (\req(pDD)) gives rise to a regular map between fiber spaces of (\req(prin)) and (\req(prin0)):
\bea(lll)
\widetilde{\varphi} : \CZ_{\widetilde{\Sigma}} \longrightarrow \CZ_{\widetilde{\Sigma}_0}, & z = (z_\delta)_{\delta \in \Lambda^{(0)} }  \mapsto \zeta = (\zeta_\upsilon)_{\upsilon \in {\cal V} ( \triangle )} & \zeta_\upsilon = \prod_{\upsilon \in F, \delta \in {\rm Int}(F) \cap \Lambda^{(0)}} z_\delta^{r^{F, \delta}_{\upsilon }} ,
\elea(zzeta)
which induces the dominating morphism $\varphi : {\PZ}_{\Sigma} \longrightarrow {\PZ}_{\Sigma_0}$  
with exceptional divisors $E_\delta$  labeled by $\delta \in \Lambda^{(0)} \setminus {\cal V} ( \triangle )$. Note the morphism $j^*: D^\dagger_{\triangle} \longrightarrow D^\dagger_{\triangle 0}$ induced by $j$ in (\req(pDD)) gives rise to the biregular morphism outside exceptional divisors: ${\PZ}_{\Sigma_0} \setminus (\cup \varphi(E_\delta) ) \stackrel{\sim}{\longrightarrow} {\PZ}_{\Sigma} \setminus (\cup E_\delta)$.

\subsection{Picard group of toric variety and anti-canonical hypersurfaces \label{ssec.Pt}} 
We shall identify an element $\nu \in {\bf n}_\triangle^\dagger$  with the character function $\chi_\nu : {\bf T} ( {\bf n}_\triangle ) \longrightarrow {\CZ}^* (= {\CZ}/{\ZZ})$. For $\nu \in {\bf n}_\triangle^\dagger$, the (orbifold) line bundle ${\cal O}(\nu)$ over ${\PZ}_{\Sigma}$ is defined through (\req(prin)) as the ${\bf T} ( {\bf n}_\triangle )$-quotient of $\CZ_{\widetilde{\Sigma}}\times \CZ$ via the action: $(p, \zeta)\cdot t = (p \cdot t, \chi_\nu (t)\zeta)$: ${\cal O}(\nu) = \CZ_{\widetilde{\Sigma}}\times_{{\bf T} ( {\bf n}_\triangle)} \CZ  $. Identifying $\nu$ with ${\cal O}(\nu)$, one may regard ${\bf n}_\triangle^\dagger$ as the Picard group ${\rm Pic}({\PZ}_{\Sigma})$ of ${\PZ}_{\Sigma}$. Indeed, we have the following result (see, e.g. (11) and Theorem 1 in \cite{R96}):
\bea(lllll)
{\bf n}^\dagger_{\triangle~\CZ} &\simeq {\rm Pic}({\PZ}_{\Sigma})_{\CZ} &\simeq H^2 ({\PZ}_{\Sigma}, {\CZ} )&, & \nu \leftrightarrow {\cal O}(\nu) \leftrightarrow {\rm Chern ~ class ~ of } ~{\cal O}(\nu). 
\elea(PicP)
An element $\varrho$ in $D^\dagger_+$ can be regarded as a divisor of $\PZ_\Sigma$, whose image under $\iota^*$ in (\req(exa*)) defines the line bundle ${\cal O}(\iota^* \varrho )$ over $\PZ_\Sigma$. The canonical bundle of $\PZ_{ \Sigma }$ is equal to ${\cal O}(\iota^* \kappa )$, where 
\be
\kappa := - \sum_{ \delta \in \Lambda^{(0)}} e^{\delta \dagger}   \in D_\triangle^\dagger 
\ele(kappa)
and $ e^{\delta \dagger}$ are basis elements in (\req(D*)). Similarly, the canonical bundle of the minimal toric variety $\PZ_{ \Sigma_0 }$ is equal to ${\cal O}(\iota_0^*  \kappa_0 )$ with 
\be
\kappa_0 =  - \sum_{\upsilon \in {\cal V} ( \triangle )} e^{\upsilon \dagger} \in D_{\triangle ~ 0}^\dagger. 
\ele(kappa0)
Under the morphism $p^*$ in (\req(p*)), $p^*(\kappa_0) = \kappa $, hence $\PZ_{ \Sigma }$ is Gorenstein and dominating  over $\PZ_{ \Sigma_0 }$. Let $X$ be a generic anti-canonical hypersurface $X$ in $\PZ_{ \Sigma }$.  Then $X$ is quasi-smooth\footnote{ Here we mean $\pi^*(X)$ is smooth hypersurface of $\CZ_{\widetilde{\Sigma}}$ in (\req(prin)), which is equivalent to the smoothness of $X$ when $n \leq 4$ (\cite{R96} Proposition 6).}, with the Picard group described by  (\cite{R96} Proposition 5\footnote{The formula (\req(PicX)) was stated in Proposition 5 of \cite{R96} for cases $n \geq 4$, but the formula is still valid for $n \leq 3$.}, or Proposition  \ref{prop:rhoRf} with $\rho= -\kappa$ in Appendix;  for the simplex-$\triangle$ case,  see \cite{R91} Lemma 7 (ii))
\be
{\rm Pic}(X)_{\CZ} \simeq {\bf n}^\dagger_{\triangle~\CZ} \oplus \bigoplus_{F, ~ \nu_F } \CZ ~ \nu_F 
\ele(PicX)
where $F$ runs over codimensional 2 faces of $\triangle$ with its dual face $F^*$, and  
$\nu_F \in ( {\rm Int}(F) \cap L ) \times ({\rm Int}(F^*) \cap L^*)$. In (\req(PicX)), the line bundles in ${\bf n}^*_{\triangle~\CZ}$  are those inherited from ${\PZ}_{\Sigma}$, and elements in the second term are the contribution from  extract exceptional divisors lying over 1-tori determined by codim-2 faces $F$ in ${\PZ}_{\Sigma_0}$ with elements in ${\rm Int}(F^*) \cap L^*$ representing primitive $0$-cohomology basis of "blow-up-centers" in the 1-torus (see,  Lemma \ref{lem:sectK} Remark (1) in Appendix). Indeed, $X$ is a pull-back of a generic hypersurface $X_0$ of $\PZ_{\Sigma_0}$, in which the closure $\overline{\rm O}_{F}$ of ${\bf T} (L)$-orbit ${\rm O}_{F}$ fixed by $L \cap \sigma (F)$ is isomorphic to $\PZ^1$.  
The line bundle ${\cal O}(\iota_0^* \kappa_0)$ over $\PZ_{\Sigma_0}$  restricting on $\overline{\rm O}_{F}$ is equivalent to ${\cal O}_{\PZ^1}(d) ~ (d = |{\rm Int}(F^*)|+1)$ over  $\PZ^1$ with the homogeneous coordinates $z=[z_1, z_2]$. Then $z^k := z_1^k z_2^{d-k} ~ (0 \leq k \leq d)$ form the monomial basis of ${\cal O}_{\PZ^1}(d)$ so that $\{ z^k \}_{k=1}^{d-1}$ are in one-to-one correspondence with ${\rm Int}(F^*) \cap L^*$. The set $X_0 \cap \overline{\rm O}_{F}$ consists of $d$ generic elements $x_j \in \PZ^1 ~ (1 \leq j \leq d)$. Let $x_j^* ~ (1 \leq j \leq d)$ be a basis of $H^0(X_0 \cap \overline{\rm O}_{F})$ dual to $x_j$'s with $\langle x_j^*, x_k \rangle = \varepsilon^j_k $. Then $p^k = \sum_{j=1}^d e^{2 \pi {\rm i}j k/d} x_j^* ~ (0 \leq k \leq d-1)$ form a basis of $H^0(X_0 \cap \overline{\rm O}_{F})$, which contains the $(d-1)$-dimensional primitive cohomology $H^{0, 0} (X_0 \cap \overline{\rm O}_{F})$ spanned by $p^k ~ (1 \leq k \leq d-1)$. Then all $\nu_F$'s in (\req(PicX)) for a fixed $(n-2)$-face $F$ are indexed by $(\delta_F, p^k )$ for $\delta_F \in {\rm Int}(F) \cap L$ and $1 \leq k \leq d-1$. 
Note that $p^0 = \sum_{j=1}^d x_j^*$, and  $(\delta_F, p^0)$ is the toric divisor of $\PZ_{\Sigma}$ determined by $\delta_F$, i.e. $e^{\delta_F \dagger} \in D^\dagger_\triangle$  in (\req(exa*)), which gives rise to the line bundle ${\cal O}(\iota^* e^{\delta_F \dagger}) \in {\bf n}^\dagger_{\triangle~\CZ}$ in (\req(PicX)). \par \vspace{.1in} \noindent
{\bf Remark.} ${\rm Pic}(X)_{\CZ}$ in (\req(PicX)) depends only on $\Lambda^{(0)}$, not on the detailed simplicial structure of $\Lambda$. The structure in Picard group which is relevant to $\Lambda$ appears only in the part ${\bf n}^\dagger_{\triangle~\CZ}$ in (\req(PicX)), reflected in the positive divisor cone discussed in the next subsection.

\subsection{Kahler cone of a toric variety\label{ssec.Kcone}} 
In this subsection, we discuss the positive cone structure in (\req(PicX)) related to ${\bf n}^\dagger_{\triangle~\QZ}$, which is a part of Kahler cone of $X$, consisting of line bundles induced by positive divisors of $\PZ_{\Sigma}$. An element in $\varrho \in D^\dagger (= D^\dagger_{\triangle~ {\QZ}})$ will be denoted by
$$
\varrho = \sum_{\delta \in \Lambda^{(0)}} \varrho^{\delta } e^{\delta \dagger} , ~ ~ \varrho^{\delta } = \langle \varrho , e^{\delta} \rangle ,
$$ 
regarded as a (toric) divisor  whose $\iota^*$-image in (\req(exa*)) defines the line bundle ${\cal O}(\iota^* \varrho)$ over $\PZ_\Sigma$ in (\req(PicP)). In the first quadrant cone of $D^\dagger$,
$$
{\rm Int}(D^\dagger_+) = \{ \varrho \in D^\dagger |  \langle \varrho , e^{\delta} \rangle > 0 ~ {\rm for} ~ \delta \in \Lambda^{(0)} \} ~ ~ \subseteq ~ D^\dagger_{+} = 
\{ \varrho \in D^\dagger  | \langle \varrho, e^{\delta} \rangle \geq 0 ~ {\rm for} ~ \delta \in \Lambda^{(0)} \} 
$$
consists of effective divisors of $\PZ_\Sigma$. First we define the following elements in ${\bf n}_{\triangle ~ \QZ}$. For $\delta \in \Lambda^{(0)}$, there is a (unique) simplex ${\tt s}_{\delta}^\vdash \in \Lambda$ so that $-\delta \in {\rm Int}({\tt s}_{\delta}^\vdash)$, equivalently, $-\delta = \sum_{\delta^\ast \in {\tt s}_{\delta}^\vdash } n^{\delta}_{\delta^\ast } \delta^\ast$ for positive rational numbers $ n^{\delta}_{\delta^\ast }$. Define 
\be 
n^\delta (= n_{\Lambda}^\delta) = e^{\delta} + \sum_{\delta^\ast \in {\tt s}_{\delta}^\vdash } n^{\delta}_{\delta^\ast } e^{\delta^\ast} \in {\bf n}_{\triangle ~ \QZ} \cap \Omega , ~ ~ (n^{\delta}_{\delta^\ast } \in \QZ_{> 0}),
\ele(ndel)
whose one-parameter group $s^{n^\delta}$ acts on ${\bf T}(D)$ in (\req(CZ)) with $(s=0)$-limit outside $\CZ_{\widetilde{\Sigma}}$. Note that the above $n^\delta$'s in (\req(ndel)) depend on the triangulation $\Lambda$, and generate a convex cone in ${\bf n}_{\triangle ~ \QZ}$, 
$$
{\bf n}_{\Lambda, 0} :=\sum_{\delta \in \Lambda^{(0)}} \QZ_{\geq 0} ~ n^\delta  \subseteq {\bf n}_{\triangle ~ \QZ} \cap \Omega.
$$
Using $n^\delta$'s in (\req(ndel)), we are going to describe a system of local coordinates of the principal bundle (\req(prin)). The following lemma is obvious:
\begin{lem}\label{lem:locS}
For a $(n-1)$-simplex ${\tt s} \in \Lambda^{(n-1)}$,  the $(d-n)$ elements, $n^{\delta'}$ for $\delta' \in  \Lambda^{(0)} \setminus {\tt s}$, form a ${\bf n}_{\triangle~\QZ}$-basis with the expression,
$$
n^{\delta'} = \sum_{\delta \in {\tt s} \cap \Lambda^{(0)}} n_\delta^{\delta'} e^\delta + \sum_{\delta'' \in \Lambda^{(0)} \setminus {\tt s}} n_{\delta''}^{ \delta'} e^{\delta''} ~ ~ {\rm with} ~  n_\delta^{\delta'},  n_{\delta''}^{ \delta'} \geq 0  ,  ~ ~ n_{\delta'}^{ \delta'} = 1 , ~ ~ ~ (\delta' \in \Lambda^{(0)}\setminus {\tt s}).
$$ 
The above ${\bf n}_{\triangle~\QZ}$-basis, together with  $e^\delta ~ (\delta \in {\tt s} \cap \Lambda^{(0)})$,  form a basis of $D_{\triangle~\QZ}$. The basis $\{ e_{\delta {\tt s} } \} \cup \{ n_{\delta', {\tt s}} \}$ of $D^\dagger$  dual to $ \{ e^\delta \} \cup \{ n^{\delta'} \}$ is related to the basis $\{ e^{\delta \dagger}\} \cup \{ e^{\delta' \dagger} \}$ in (\req(D*))  by
\bea(ll)
e_{\delta , {\tt s}} = e^{\delta \dagger} + \sum_{\delta'' \in \Lambda^{(0)} \setminus {\tt s}} g_{\delta }^{\delta''} e^{\delta'' \dagger}, & 
n_{\delta', {\tt s}} = \sum_{\delta'' \in \Lambda^{(0)} \setminus {\tt s}} g_{\delta' }^{\delta''} e^{\delta'' \dagger}  ; \\
e^{\delta \dagger} = e_{\delta , {\tt s}}  + \sum_{\delta' \in \Lambda^{(0)} \setminus {\tt s}} n_{\delta }^{\delta'} n_{\delta' {\tt s}} , &
e^{\delta'' \dagger} = \sum_{\delta' \in \Lambda^{(0)} \setminus {\tt s}} n_{\delta'' }^{\delta'} n_{\delta' {\tt s}}, 
\elea(duals)
where the index $\delta \in {\tt s} \cap \Lambda^{(0)}$, $
\delta', \delta'' \in \Lambda^{(0)} \setminus {\tt s}$, and the coefficients satisfy the relations, $(g_{\delta }^{\delta''} ) = - (n_{\delta }^{\delta'})(n_{\delta' }^{\delta''})^{-1}$ and $( g_{\delta' }^{\delta''}) = (n_{\delta' }^{\delta''})^{-1}$. Furthermore, $e_{\delta, {\tt s}} ~ (\delta \in {\tt s} \cap \Lambda^{(0)})$ form a rational basis of the subspace $\beta^*(L^*_{\QZ})$ of $D^\dagger$ in (\req(exa*)).
\end{lem} 
\par \vspace{.1in} \noindent
For a $(n-1)$-simplex ${\tt s} \in \Lambda^{(n-1)}$, the generators ${\tt s} \cap \Lambda^{(0)}$ of the simplicial cone $\sigma ({\tt s})$ span a sublattice $L_{\tt s}$ of $L$. It is known that the collection of  $U_{\tt s} = {\rm Spec}[ \widehat{\sigma ({\tt s})} \cap L^* ]$ for  ${\tt s} \in \Lambda^{(n-1)}$ forms an affine open chart of ${\PZ}_{\Sigma}$.  For a given ${\tt s}$, the ${\bf T} ( {\bf n}_\triangle )$-space  (\req(prin)) over $U_s$ is given by $\pi^{-1}(U_{\tt s}) = {\rm Spec}[ \widehat{\widetilde{\sigma} ({\tt s})} \cap D_\triangle^\dagger ]$. Consider the basis $\{ e^\delta \} \cup \{ n^{\delta'} \}$ of $D_{\triangle~\QZ}$ in Lemma \ref{lem:locS}. Then  $e^\delta$'s and $n^{\delta'}$'s generate a cone $\Omega_{\tt s}$ in $D_{\triangle~\RZ}$ which contains $\widetilde{\sigma} ({\tt s})$ as a $n$-face. The relations of cones and lattices, $(\widetilde{\sigma} ({\tt s}), D_\triangle) \hookrightarrow (\Omega_{\tt s}, D_\triangle) \stackrel{\beta}{\longrightarrow} (\sigma ({\tt s}), L)$, give rise to the morphisms of toric varieties: $\pi^{-1}(U_{\tt s}) \hookrightarrow \overline{\pi^{-1}(U_{\tt s})} \stackrel{\overline{\pi}}{\longrightarrow} U_{\tt s}$, where $\beta$ is the morphism in (\req(exa)), and $\overline{\pi^{-1}(U_{\tt s})}= {\rm Spec}[ \widehat{\Omega_{\tt s}} \cap D_\triangle^\dagger ]$. For $\delta' \notin {\tt s}$, let $k_{\delta'}$ be the positive integer characterized by the primitive property of $k_{\delta'} n^{\delta'}$ in ${\bf n}_\triangle$. Consider the lattices , 
$$
\begin{array}{ll}
{\bf n}_{\triangle, {\tt s}} := \sum_{\delta' \in \Lambda^{(0)}\setminus {\tt s}} \ZZ ( k_{\delta'} n^{\delta'}) \subseteq
{\bf n}_{\triangle}, & D_{\triangle, {\tt s}}: = \sum_{\delta \in {\tt s} \cap \Lambda^{(0)}} \ZZ e^{\delta} + {\bf n}_{\triangle, {\tt s}} \subseteq D_\triangle.
\end{array}
$$
The cone-lattice relation $(\Omega_{\tt s}, D_{\triangle, {\tt s}}) \stackrel{\beta}{\longrightarrow} (\sigma ({\tt s}), L_{\tt s})$ induces the toric morphism, $\CZ^d = \CZ^n \times \CZ^{d-n} \longrightarrow \CZ^n$, with the coordinates corresponding to $e_{\delta , {\tt s}}$'s and $\frac{n_{\delta', {\tt s}}}{k_{\delta'}}$'s in Lemma \ref{lem:locS}.  The morphism $\overline{\pi}$ is equivalent to the projection from $\CZ^d / (D_\triangle/ D_{\triangle, {\tt s}} ) (= \overline{\pi^{-1}(U_{\tt s})})$ to $\CZ^n / (L/ L_{\tt s}) (= U_{\tt s})$, and  we have $\pi^{-1}(U_{\tt s}) = (\CZ^n \times \CZ^{* n-d})  / (D_\triangle/ D_{\triangle, {\tt s}} )$ and ${\bf T}({\bf n}_\triangle) = \CZ^{* n-d}/({\bf n}_\triangle/{\bf n}_{\triangle, {\tt s}})$. Using $e^{\delta} ~ (\delta \in {\tt s} \cap \Lambda^{(0)})$, one can lift $L$ into $D_\triangle$ so that $D_\triangle \simeq L \oplus {\bf n}_\triangle.$ Hence the inclusion $\pi^{-1}(U_{\tt s}) \hookrightarrow \overline{\pi^{-1}(U_{\tt s})}$ is equivalent to
$$
\begin{array}{llll}
 \pi^{-1}(U_{\tt s}) \simeq U_{\tt s} \times {\bf T}({\bf n}_\triangle) & \hookrightarrow &U_{\tt s} \times \overline{{\bf T}({\bf n}_\triangle)} \simeq \overline{\pi^{-1}(U_{\tt s})}, & \overline{{\bf T}({\bf n}_\triangle)} := \CZ^{n-d}/({\bf n}_\triangle/{\bf n}_{\triangle, {\tt s}} ) ,
\end{array}
$$
and $\pi$ in (\req(prin)) corresponds to the projection to the first component.

For a simplex ${\tt s} \in \Lambda^{(n-m_{\tt s}-1)}$ for $0 \leq m_{\tt s} \leq n-1$, it is known that the ${\bf T} (L)$-orbit ${\rm O}_{\tt s}$ fixed by $ L \cap \sigma ({\tt s})$ is isomorphic to a $m_{\tt s}$-torus, whose closure $\overline{\rm O}_{\tt s}$  is a $m_{\tt s}$-dimensional toric variety \cite{Dan, KKMS, Oda} (or in Appendix (\req(Os))  for the detailed structure). Consider the subspace of $D^\dagger$,
$$
D^\dagger_{\tt s} = \{ \varrho \in D^\dagger  | \langle \varrho , e^{\delta} \rangle = 0  ~ 
{\rm for} ~ \delta \in {\tt s} \},
$$
in which, there is the positive quadrant cone:
$$
{\rm Int}(D^\dagger_{{\tt s}, +})  = \{ \varrho \in D^\dagger_{\tt s}  | \langle \varrho , e^{\delta'} \rangle > 0 ~ {\rm for} ~ \delta' \notin {\tt s} \} ~ ~ \subseteq ~ D^\dagger_{{\tt s}, +} = 
\{ \varrho \in D^\dagger_{\tt s}  | \langle \varrho , e^{\delta'} \rangle \geq 0 ~ {\rm for} ~ \delta' \notin  {\tt s} \}.
$$
Then  $\pi^{-1}(\overline{\rm O}_{\tt s})$ in the fibration (\req(prin)) is contained in $D^\dagger_{{\tt s} ~ \CZ}$ with the coordinate system given by $\{ e^{\delta' \dagger} \}_{\delta' \notin {\tt s}} \subset 
D^\dagger_{{\tt s}, +}$.  For a $(n-1)$-simplex ${\tt s} \in \Lambda^{(n-1)}$, using the local system $\{ e_{\delta , {\tt s} } \} \cup \{ n_{\delta', {\tt s}} \}$ in Lemma \ref{lem:locS}, one can express an element $\varrho$ in $D^\dagger$ 
by 
\bea(lll)
\varrho =  \varrho_{\tt s}^\prime + \varrho_{\tt s}, &\varrho_{\tt s}^\prime:= \sum_{\delta \in {\tt s} \cap \Lambda^{(0)}} \langle \varrho , e^{\delta} \rangle e_{\delta , {\tt s}} \in \beta^*(L^*_{\QZ}) , &\varrho_{\tt s} := \sum_{\delta'  \in \Lambda^{(0)} \setminus {\tt s} }\langle \varrho , n^{\delta'} \rangle n_{\delta', {\tt s}} \in D^\dagger_{\tt s}.
\elea(rhos)
In the above expression, $\varrho_{\tt s}^\prime$ can be regarded as the linear functional of $L_Q$ which takes the value $\rho^\delta $ on $\delta \in {\tt s} \cap \Lambda^{(0)}$, hence $\langle \varrho_{\tt s}, e^{\delta''} \rangle$ is the difference of $e^{\delta''} $-values between $\varrho$ and  $L_{\rm \QZ}$-linear functional $\varrho_{\tt s}^\prime$  for  $\delta'' \notin {\tt s}$.  \par \noindent
{\bf Definition.}\cite{O88, OP, R96}\footnote{Here the definition is defined via the divisor lattice $D_{\triangle, ~ \QZ}$, which is equivalent to the notion in \cite{O88, OP, R96} where the convex divisor is formulated through the graph over $L_{\QZ}$. Indeed, the $l_\sigma, f_\sigma$ in Definition 3 of \cite{R96} are corresponding to $\varrho_{\tt s}^\prime, \varrho$ here with  $\sigma = \sigma ({\tt s})$. Furthermore, results and definitions of Subsection  \ref{ssec.Kcone} here can be carried over to a general toric variety $\PZ_{\Sigma}$ defined by a simplicial cone decomposition of $L_{\QZ}$. Then $\Lambda^{(0)}$ is identified with the collection of all primitive $L$-elements representing  1-dimensional cones in $\Sigma$, but the corresponding  $\triangle$ in $L_{\RZ}$ may not be convex in general.}   Let $\varrho$ be an element in $D^\dagger (= D^\dagger_{\triangle~ {\QZ}})$. \par \noindent
(I) $\varrho$ is a non-negative (or convex) divisor if $\varrho_{\tt s}$ in (\req(rhos)) is an element of $D^\dagger_{{\tt s}, +}$ for all ${\tt s} \in \Lambda^{(n-1)}$. \par \noindent
(II) $\varrho$ is a positive (or strictly convex) divisor if $\varrho_{\tt s} \in {\rm Int}(D^\dagger_{{\tt s}, +})$ for all ${\tt s} \in \Lambda^{(n-1)}$.
\par \vspace{.1in} \indent
Consider the positive cone in $D^\dagger$, 
\be
C^\dagger_+ ~ (= C^\dagger_{\Lambda, +}) = \{ \varrho \in  D^\dagger | \langle \varrho_{\tt s}, e^{\delta'} \rangle \geq 0 ~ ~ {\rm for} ~ \delta' \in \Lambda^{(0)} \setminus {\tt s}, ~ {\tt s} \in \Lambda^{(n-1)}  \},
\ele(Pcone)
whose interior consists of all positive divisors whose $\iota^*$-images in (\req(exa*)) are positive line bundles over $\PZ_\Sigma$. We now derive an mechanism to express the positive cone $C^\dagger_+$. For a $(n-1)$-simplex ${\tt s} \in \Lambda^{(n-1)}$ and $\delta' \in \Lambda^{(0)} \setminus {\tt s}$, we consider the following element in ${\bf n}_{\triangle ~ \QZ}$:
\be 
n^{\delta'}_{\tt s} := e^{\delta'} + \sum_{\delta \in {\tt s} \cap \Lambda^{(0)}} n^{\delta'}_{\delta, {\tt s}}e^{\delta} ~ ~ ~ ~ {\rm with} ~ ~ \delta' + \sum_{\delta \in {\tt s} \cap \Lambda^{(0)}} n^{\delta'}_{\delta, {\tt s}}\delta = 0 \in L_{\rm \QZ}.
\ele(nsdel)
As ${\tt s} \cap \Lambda^{(0)}$ form a basis of $L_{\QZ}$, the above element is uniquely defined. Note that $n^{\delta'}_{\tt s}= n^{\delta'}$ in (\req(ndel)) iff ${\tt s}_{\delta'}^\vdash \subseteq {\tt s}$. Furthermore, 
$ n^{\delta'_2}_{{\tt s}_1}$ is propositional to $n^{\delta'_1}_{{\tt s}_2}$ by a positive number if ${\tt s}_1 \cap {\tt s}_2$ is a common codim-1 face of ${\tt s}_i$ , and  $\delta'_i \in {\tt s}_i$ for $i=1, 2$. 
Define the positive cone in ${\bf n}_{\triangle ~ \QZ}$,
\be 
{\bf n}_{\Lambda, +} :=\sum_{{\tt s} \in \Lambda^{(n-1)}, \delta' \in \Lambda^{(0)} \setminus {\tt s}} \QZ_{\geq 0} ~ n^{\delta'}_{\tt s}  \subseteq {\bf n}_{\triangle ~ \QZ}. 
\ele(n+)
Consider the dual of ${\bf n}_{\Lambda, +}$,
\be
{\bf n}^\dagger_{\Lambda, +}= {\rm dual ~ cone ~ of ~ } {\bf n}_{\Lambda, +}  ~ {\rm in} ~  {\bf n}^\dagger_{\triangle~ {\rm \QZ}} , ~ ~ ~ ~
\widehat{\bf n}_{\Lambda, +} (D^\dagger)={\rm dual ~ cone ~ of ~} {\bf n}_{\Lambda, +} ~ {\rm in} ~  D^\dagger_{\triangle~{\rm \QZ}}. 
\ele(ndual)
Through $\iota^*$ in (\req(exa*)), the cones in (\req(ndual)) are related by 
$$
\widehat{\bf n}_{\Lambda, +} (D^\dagger) = \iota^{* ~ -1}({\bf n}^\dagger_{\Lambda, +}), ~ ~ {\rm Int}(\widehat{\bf n}_{\Lambda, +} (D^\dagger)) = \iota^{* ~ -1}({\rm Int}({\bf n}^\dagger_{\Lambda, +})).
$$ 
\begin{prop}\label{prop:+cone} 
{\rm (I)} The positive cone $C^\dagger_+$ in (\req(Pcone)) is equal to $\widehat{\bf n}_{\Lambda, +} (D^\dagger)$. \par \noindent
{\rm (II)} The non-negative divisors in the first quadrant cone of $D^\dagger$ are given by
$$
\begin{array}{l}
C^\dagger_+ \cap D^\dagger_{+} = \{ \varrho \in D^\dagger_{+} | ~ \langle \varrho,  n^{\delta'}_{\tt s} \rangle \geq 0 ~ 
~ {\rm for} ~ \delta' \notin  {\tt s} \in \Lambda^{(n-1)},  ~ {\tt s}_{\delta'}^\vdash \setminus {\tt s} \neq \emptyset \},\\
{\rm Int}(C^\dagger_+) \cap {\rm Int}(D^\dagger_+) = \{ \varrho \in {\rm Int}(D^\dagger_+) | ~ \langle \varrho,  n^{\delta'}_{\tt s} \rangle > 0 ~ 
~ {\rm for} ~ \delta' \notin  {\tt s} \in \Lambda^{(n-1)},  ~ {\tt s}_{\delta'}^\vdash \setminus {\tt s} \neq \emptyset \} , 
\end{array}
$$
which satisfy the relations,  $\iota^*( C^\dagger_+ \cap D^\dagger_{+}) = {\bf n}^\dagger_{\Lambda, +}$ and 
$\iota^*( {\rm Int}(C^\dagger_+) \cap {\rm Int}(D^\dagger_+)) = {\rm Int} ({\bf n}^\dagger_{\Lambda, +})$.
\end{prop}
{\it Proof.} For an element $\delta' \notin  {\tt s} \in \Lambda^{(n-1)}$, using Lemma \ref{lem:locS} and (\req(rhos)), one finds 
$$
\langle \varrho_{\tt s}, e^{\delta'} \rangle = \sum_{\delta \in {\tt s}} (  \sum_{\delta''  \notin  {\tt s} }  n^{\delta''}_{\delta} g_{\delta''}^{\delta'})\varrho^{\delta} + \varrho^{\delta'}= \langle \varrho,  \sum_{\delta \in {\tt s}} (\sum_{\delta''  \notin  {\tt s} }  n^{\delta''}_{\delta} g_{\delta''}^{\delta'})e^{\delta} + e^{\delta'}\rangle  = \langle \varrho, \sum_{\delta'' \notin {\tt s}} n^{\delta''}g^{\delta'}_{\delta''} \rangle  = \langle \varrho, n_{\tt s}^{\delta'} \rangle .
$$
The last equality in above follows from the form of $n^{\delta''} \in {\bf n}_{\triangle ~{\rm \QZ}}$ and the definition of $n_{\tt s}^{\delta'}$ in (\req(nsdel)). Hence we obtain (I). When ${\tt s}_{\delta'}^\vdash \subseteq {\tt s}$, $n^{\delta'}_{\tt s}= n^{\delta'}$ with $n^{\delta'}_{\delta, {\tt s}}$'s in (\req(nsdel)) all non-negative, which implies the constraint $\langle \varrho_{\tt s}, e^{\delta'} \rangle \geq 0$ (or $>0$) for $C^\dagger_+ \cap D^\dagger_{+}$ ( or ${\rm Int}(C^\dagger_+) \cap {\rm Int}(D^\dagger_+ )$ respectively) is redundant. Since $\beta^*(L^*_{\rm \QZ})$ is contained in $C^\dagger_+ $ with $n_{\tt s}^{\delta'}$-zero value,  the relation between divisors and line bundles in (II) can be derived from the definition of non-negative or positive divisors.
$\Box$ \par \vspace{.2in} \noindent
By Proposition \ref{prop:+cone}, ${\bf n}^\dagger_{\Lambda, +}$ is Kahler cone of $\PZ_{\Sigma}$ consisting of all non-negative line bundles over $\PZ_{\Sigma}$, regarded as  a part of Kahler cone of $X$, and will be denoted by
\be
{\cal C}_{\rm Pic} (X) (= {\cal C}_{{\rm Pic}, \Lambda}(X) )= {\bf n}^\dagger_{\Lambda, +}.
\ele(Kcone)
The vertices of ${\bf n}^\dagger_{\Lambda, +}$, dual to facets of ${\bf n}_{\Lambda, +}$ in (\req(n+)), can be described by a lifting of ${\bf n}^\dagger_{\Lambda, +}$ in $C^\dagger_+$. We now provide some examples as demonstration of the computation of ${\bf n}^\dagger_{\Lambda, +}$; the first two are well-known, which serve as useful 'toy'-models to illustrate the method and notations used in this subsection.

{\bf Example 1.}\footnote{In this example, the conditions, (\req(Lam0)) and (\req(L0)),  are not satisfied when $G \neq 1$. The positive divisors discussed here are those for the minimal toric variety $\PZ_{(\Sigma_0, L)}$. }  
${\rm Rank}({\bf n}_\triangle)= 1$ case, i.e. $(\triangle, L)$ = a reflexive simplex. Equivalently, $\PZ_{\Sigma_0}$  is isomorphic to an abelian quotient of weighted projective space \cite{B, R94, R96}\footnote{Some misprints occur in RIMS journal paper \cite{R96} page 832, where $d_j:= \frac{d}{n_j}$ and ${\rm dia.}[ e^{2\pi i/d_1}, \ldots , e^{2\pi i/ d_5}]$ (on line 14, 15) should be $d_j:= \frac{d}{{\rm g.c.d}(d, n_j)}$, ${\rm dia.}[ e^{2\pi i n_1/d}, \ldots , e^{2\pi i n_5/ d}]$ respectively.}: $\PZ_{\Sigma_0} \simeq {\PZ}_{(w_i)}^{n}/G$, with weights $w_i \in \ZZ_{>0} ~ (1\leq i \leq n+1)$, and $G$ a subgroup of the diagonal finite group $SD ~ (:=\{{\rm dia}[ t_1, \ldots , t_{n+1}] | t_i^{d_i} = \prod_{i}t_i=1  ~ ( 1 \leq i \leq n+1)\})$, containing the element ${\rm dia}
[ e^\frac{2\pi {\rm i} w_1}{d}, \ldots , e^\frac{2\pi {\rm i}w_{n+1}}{d}] \in G$, where $d= \sum_{i=1}^{n+1} w_i$ and $d_i = \frac{d}{{\rm g.c.d}(d, w_i)}$. Indeed, the simplex $\triangle$ is spanned by  $(n+1)$ generators $\delta^{i}$'s of a $n$-dimensional real vector space $V (\simeq \RZ^n)$ satisfying the relation $\sum_{i=1}^{n+1} w_i \delta^i = 0$. Consider the $n$-lattices in $V$:
$L_0 = \sum_{i=1}^{n+1} \ZZ \delta^{i}$, and $\widehat{L} = \{ \sum_{i=1}^{n+1} \frac{k_i}{d_i}\delta^{i} | ~ k_i,  \sum_i \frac{k_i}{d_i}  \in \ZZ ~ \forall i \}$.  Then $L$ is the lattice in $V$ satisfying the relations, $L_0 \subseteq L \subseteq \widehat{L}$ and $L/L_0 \simeq G$. The sublattice ${\bf n}_{\triangle, L_0}$ of $D_{\triangle, L_0}$ in (\req(exa)) is given by $\ZZ {\tt n} $ with ${\tt n}= \sum_{i=1}^{n+1} w_i e^{\delta^i}$, and homogeneous coordinates in (\req(prin0)) are given by $\zeta =(\zeta_i)_{i=1}^{n+1}$  with $\zeta_i$ corresponding to $e^{\delta^{i \dagger} } \in D_{\triangle, L_0}^\dagger$. The triangulation $\Lambda_0$ is the face-decomposition of $\triangle$ with the $(n-1)$-simplex ${\tt s}_i$ spanned by $\{ \delta^j \}_{1 \leq j \neq i \leq n+1}$. Then $n^{\delta^i}_{{\tt s}_i} = {\tt n}$ for all $i$, hence ${\bf n}_+ = \QZ_{\geq 0} {\tt n}$, and ${\bf n}^\dagger_+ = \QZ_{\geq 0} {\tt n}^\dagger $, where ${\tt n}^\dagger$ is the ${\bf n}_\triangle^\dagger$-basis dual to  ${\tt n}$ with $\langle {\tt n}^\dagger, {\tt n} \rangle = 1$.  The element ${\tt n}^\dagger$ gives rise to a line bundle over $\PZ_{\Sigma_0}$, denoted by ${\cal O}(1)$. The positive cone $C^\dagger_+$ consists of $\varrho \in D^\dagger$ with $\sum_i w_i \varrho^i \geq 0$, whose divisor class is the non-negative line bundle  ${\cal O}(\sum_i w_i)$ over $\PZ_{\Sigma_0}$, with the relation $C^\dagger_+ \cap {\rm Int}(D^\dagger_+) =  {\rm Int}(D^\dagger_+)$. In particular, the divisor defined by $\zeta_i=0$ gives rise to the bundle ${\cal O}(w_i)$.    

{\bf Example 2.}\footnote{The polytope in this example is not reflexive. As stated in footnote of the definition of positive divisors, Proposition \ref{prop:+cone} also holds for a general toric variety.} Some ${\rm rank}({\bf n}_\triangle)= 2$ case. We consider the blow-up of weighted projective space ${\PZ}_{(w_i)}^{n}$ at the point $c^k =(\varepsilon_i^k)_{i=1}^{n+1}$, say $k=n+1$, denoted by 
$$
\pi: \widehat{\PZ}_{(w_i)}^{n} \longrightarrow {\PZ}_{(w_i)}^{n}, ~ ~ \pi^{-1}(c^{n+1}) = {\rm exceptional ~ divisor} ~ E .
$$ 
With $\delta^{i} ~ (1 \leq i \leq n+1)$, $L = L_0 $ in Example 1, and $\delta^{n+2}:= - \delta^{n+1}$, the polytope $\triangle$ is the convex hull spanned by $\{ \delta^{i} \}_{i=1}^{n+2}$, and $\Lambda$ is simplicial decomposition of $\partial \triangle$ so that $\Lambda^{(n-1)}$ consists of the first $n$ simplices ${\tt s}_i ~ (1 \leq i\leq n) $ in Example 1, together with $n$ simplices ${\tt t}_i ~ (1 \leq i \leq n)$ defined by  ${\tt t}_i = $ the $(n-1)$-simplex spanned by $\delta^k ~ ( 1 \leq k \neq i \leq n)$ and $\delta^{n+2}$. Then $\PZ_{\Sigma} \simeq \widehat{\PZ}_{(w_i)}^{n}$. Since $\sum_{i=1}^{n+1} w_i \delta^i = 0$ and $\delta^{n+1}+ \delta^{n+2} = 0$, the subspace ${\bf n}_{\triangle ~ \QZ} $ of $D_{\triangle ~ \QZ}$ in (\req(exa)) is generated  by  ${\tt n}= \sum_{i=1}^{n+1} w_i e^{\delta^i}$ and ${\tt n}^1= e^{\delta^{n+1}}+ e^{\delta^{n+2}}$.  The vectors in (\req(nsdel)) are
$$
n^{\delta^i}_{{\tt s}_i} = {\tt n} , ~ ~ n^{\delta^{n+2}}_{{\tt s}_i} = n^{\delta^{n+1}}_{{\tt t}_i} = {\tt n}^1 , ~ ~
n^{\delta^i}_{{\tt t}_i} = {\tt n}^2, ~ ~   (1 \leq i \leq n) 
$$
where ${\tt n}^2 = {\tt n} - w_{n+1} {\tt n}^1$.  The positive cones in (\req(n+)), (\req(ndual)) and Proposition \ref{prop:+cone} are now given by 
$$
\begin{array}{c}
{\bf n}_+ =  \QZ_{\geq 0} {\tt n}^1 + \QZ_{\geq 0} {\tt n}^2 , ~ ~ ~ ~ ~ {\bf n}^\dagger_+ = \QZ_{\geq 0} {\tt n}_1 + \QZ_{\geq 0} {\tt n}_2, \\
C^\dagger_+ = \{ \varrho \in D^\dagger | \varrho^{n+1} + \varrho^{n+2} \geq 0 , \sum_{i=1}^{n-1} w_i \varrho^i -
 w_{n+1}  \varrho^{n+2} \geq 0 \}, \\

{\rm Int}(C^\dagger_+) \cap {\rm Int}(D^\dagger_+) = \{ \varrho \in {\rm Int}(D^\dagger_+) |  \sum_{i=1}^{n-1} w_i \varrho^i - w_{n+1}  \varrho^{n+2} > 0 \}, 
\end{array}
$$
where ${\tt n}_1 , {\tt n}_2$ are the basis of ${\bf n}^\dagger_{\triangle ~ \QZ} $ dual to ${\tt n}^1, {\tt n}^2$. The morphism $\iota^*$ in (\req(exa*)) is given by $
\varrho  \mapsto (\varrho^{n+1} + \varrho^{n+2}){\tt n}_1 +  (\sum_{i=1}^{n-1} w_i \varrho^i -w_{n+2}\varrho^{n+2}) {\tt n}_2$,
which sends $C^\dagger_+$ to ${\bf n}^\dagger_+$. Since $\iota^*(e^{\delta^{n+1}  \dagger} )= {\tt n}_1$ and $\iota^*(e^{\delta^1 \dagger})= w_1 {\tt n}_2$, we find  the non-negative line bundles over $\widehat{\PZ}_{(w_i)}^{n}$:
${\tt n}_1 = \pi^*{\cal O}(w_{n+1})$ trivial on the exceptional divisor $E$, and ${\tt n}_2 = \pi^*{\cal O}(w_1) \otimes {\cal O}(-w_1 E)$ trivial on fibers of the canonical projection from $\widehat{\PZ}_{(w_i)}^{n}$ to ${\PZ}_{(w_i | 1 \leq i \leq n)}^{n-1}$.

{\bf Example 3.} Flop of hypersurfaces in 4-dimensional toric variety. Let $\PZ_{\Sigma(\Lambda)}$ be a 4-dimensional toric variety associated to a triangulation $\Lambda$ of $\partial \triangle$ satisfying (\req(Lam0)) and (\req(L0)). Suppose that there are two 3-simplices ${\tt s}, {\tt t} \in \Lambda^{(3)}$ with the form ${\tt s} = \overline{\delta^0 \delta^1 \delta^2 \delta^3}, ~ {\tt t} = \overline{\delta^0 \delta^4 \delta^2 \delta^3 }$, where $\delta^i ~ ( 1 \leq 4)$  are 4 elements in a 2-dimensional face of $\triangle$ satisfying the relation $\delta^1+ \delta^4 = \delta^2 + \delta^3$. Then a generic anti-canonical hypersurface $X$ of $\PZ_{\Sigma(\Lambda)}$ is non-singular in the open chart $U_{\tt s} \cup U_{\tt t}$ of $\PZ_{\Sigma(\Lambda)}$. Define another two 3-simplices, ${\tt s}^\natural= \overline{\delta^0\delta^1 \delta^2\delta^4}$ and ${\tt t}^\natural= \overline{\delta^0\delta^1\delta^3\delta^4}$, in $\partial \triangle$ with the property ${\tt s} \cup {\tt t} = {\tt s}^\natural \cup {\tt t}^\natural$.  We form another triangulation $\Lambda^\natural$ of $\partial \triangle$ by replacing ${\tt s}, {\tt t} \in \Lambda^{(3)}$ for ${\tt s}^\natural, {\tt t}^\natural \in \Lambda^{\natural ~(3)}$. The toric variety $\PZ_{\Sigma(\Lambda^\natural)}$ is a flop of $\PZ_{\Sigma(\Lambda)}$ so that the anti-canonical hypersurface $X^\natural$ of $\PZ_{\Sigma(\Lambda^\natural)}$ defined by the same equation as $X$ is a flop of $X$ in birational geometry of 3-folds. In $D^\dagger$, there are two positive cones, $C^\dagger_{\Lambda, +}$ for $\PZ_{\Sigma(\Lambda)}$ and $C^\dagger_{\Lambda^\natural, +}$ for $\PZ_{\Sigma(\Lambda^\natural)}$. Since $\delta^1+ \delta^4 = \delta^2 + \delta^3$,  ${\tt n}:= e^{\delta^1}+ e^{\delta^4} - e^{\delta^2} - e^{\delta^3}$ is an element in ${\bf n}_\triangle$, and  vectors in (\req(nsdel)) are expressed by
$$
\begin{array}{ll}
\Lambda:  n_{\tt s}^{\delta^4} = n_{\tt t}^{\delta^1}= {\tt n}, &
\Lambda^\natural :   n_{{\tt s}^\natural}^{\delta_3} = n_{{\tt 4}^\natural}^{\delta_2} = - {\tt n}. 
\end{array}
$$
By Proposition \ref{prop:+cone}, we find
$$
\begin{array}{ll}
C^\dagger_{\Lambda, +} \subset \{ \varrho \in D^\dagger | \langle \varrho, {\tt n} \rangle \geq 0 \} , &
C^\dagger_{\Lambda^\natural, +} \subset \{ \varrho \in D^\dagger | \langle \varrho, {\tt n} \rangle \leq 0 \},
\end{array}
$$
hence $\PZ_{\Sigma(\Lambda)}$ and $\PZ_{\Sigma(\Lambda^\natural)}$ have no common positive divisors, equivalently ${\rm Int}({\bf n}^\dagger_{\Lambda, +}) \cap {\rm Int}({\bf n}^\dagger_{\Lambda^\natural, +}) = \emptyset$.
The process of flops indeed appears  in the Gorenstein toric variety $\PZ_{\Sigma(\Lambda)}= \widehat{\PZ}_{(w_i)}^{4}/G$ over ${\PZ}_{(w_i)}^{4}/G$ with $\Lambda$ satisfying (\req(Lam0)) and (\req(L0)). Since any triangulation $\Upsilon$ on all faces of codimension $\leq 2$ of a $n$-simplex $\triangle$  satisfying (\req(Lam0)) can always be extended to a triangulation $\Lambda$ of $\partial \triangle$ with the same $\Lambda^{(0)} = \Upsilon^{(0)}$, one can construct a flop process by  changing the triangulation on 2-dimensional faces of the 4-dimensional simplex $\triangle$ whenever the group $G$ is sufficient "big" , e.g. $\widehat{\PZ}^{4}/SD$ with $G= SD$ and $w_i=1$ for all $i$. Indeed, the extended $(n-1)$-simplices $\Lambda^{(n-1)}$ of $\partial \triangle$ can be constructed from $\Upsilon$ as follows. Let $\triangle_{n-1}$ be a facet of $\triangle$. If $\triangle_{n-1} \cap \Lambda^{(0)}$ is equal to ${\cal V} (\triangle_{n-1})$ (the set of vertices of $\triangle_{n-1}$),  we set ${\tt s} = \triangle_{n-1} \in \Lambda^{(n-1)}$. Otherwise, $\triangle_{n-1}$ can be expressed as the join of a vertex $v$ and a $(n-2)$-face  $\triangle_{n-2}$ of $\triangle_{n-1}$, $\triangle_{n-1} = v \ast  \triangle_{n-2}$,  with the property $\Lambda^{(0)} \setminus {\cal V} (\triangle_{n-2}) \neq \emptyset$.  
Write $\triangle_{n-2} = \cup {\tt u}$, where ${\tt u} \in \Upsilon^{(n-2)}$ and ${\tt u} \subset \triangle_{n-2}$. 
Then $\triangle_{n-1} = \cup (v \ast {\tt u})$. We set ${\tt s}= v \ast {\tt u} \in \Lambda^{(n-1)}$ if ${\tt u} \subset {\rm Int}(\triangle_{n-2})$. Otherwise, 
${\tt u}$ is spanned by a finite number of vertices of $\triangle_{n-2}$, $v_1, \cdots, v_{n-1-k} ~ (k \geq 1)$, together with $\delta_1, \cdots, \delta_k \in \Lambda^{(0)} \setminus {\cal V} (\triangle_{n-2})$.  Then $
v \ast {\tt u} = \triangle_{n-1-k}^\prime \ast \overline{\delta_1 \ldots, \delta_k}$, where $\triangle_{n-1-k}^\prime$ is the $(n-1-k)$-face of $\triangle$ spanned by $v, v_1, \ldots , n_{n-1-k}$. One finds $\triangle_{n-1-k}^\prime = \cup {\tt u}^\prime$ with ${\tt u}^\prime \in \Upsilon^{(n-1-k)}$, and $v \ast {\tt u} = \cup ({\tt u}^\prime \ast \overline{\delta_1 \ldots, \delta_k})$, by which we set ${\tt s}= {\tt u}^\prime \ast \overline{\delta_1 \ldots, \delta_k} \in \Lambda^{(n-1)}$. Therefore we obtain a triangulation $\Lambda$ of $\partial \triangle$ as an extension of $\Upsilon$  with the same 0-skeleton.

\section{Deformation Space of Anti-canonical Hypersurface in Toric Variety\label{sec.Defm}}
\setcounter{equation}{0}
Consider the dual polytope $(\triangle^*, L^*)$ in $L^*_{\RZ}$, equipped with a simplicial decomposition $\Lambda^* = \{ \Lambda^{* (j)} \}_{j=1}^{n-1}$ of $\partial \triangle^*$ satisfying (\req(Lam0)) and (\req(L0)) for $\Lambda^*$ and $L^*$:
\bea(lll)
0 \longrightarrow {\bf n}_{\triangle^*} \stackrel{\iota'}{\longrightarrow} D_{\triangle^*} \stackrel{\beta'}{\longrightarrow} L^* \longrightarrow 0, & ~ &
0 \longrightarrow L \stackrel{\beta^{\prime *}}{\longrightarrow} D_{\triangle^*}^\dagger \stackrel{\iota^{\prime *}}{\longrightarrow} {\bf n}_{\triangle^*}^\dagger \longrightarrow  0  , 
\elea(exdual)
where $ \beta'( e^{\delta^*}) = \delta^*$ for $\delta^* \in \Lambda^{* (0)}$. Note that no relation is imposed between the triangulation $\Lambda^*$ of $\partial \triangle^*$ and $\Lambda$ of $\partial \triangle$ in Section \ref{sec.PG}. We shall denote the origin of $L^*_{\RZ}$ by $0^*$, and the collection of vertices of $\triangle^*$ by ${\cal V} ( \triangle^* )$. A simplex in $\Lambda^*$ will be denoted by ${\tt s}^*$, and $d^*:= | \Lambda^{(0)} |$.  Corresponding to the basis $e^{\delta^*}$'s of $D_{\triangle^*}$, the torus ${\bf T}(D_{\triangle^*}) (\simeq \CZ^{*  d^*})$ is parametrized by coordinates $u = (u_{\delta^*})_{\delta^* \in \Lambda^{* (0)}}$. The simplicial decomposition $\Lambda^*$ of $\partial \triangle^*$ gives rise to the ${\bf T}(D_{\triangle^*})$-variety $\CZ_{\widetilde{\Sigma^*}}$ as in (\req(CZ)):
\be
u = (u_{\delta^*})_{\delta^* \in \Lambda^{* (0)}} \in  \CZ_{\widetilde{\Sigma^*} } = D_{\triangle^*~\CZ} -  \bigcup_{K \notin \{ {\tt s} \cap \Lambda^{* (0)} | {\tt s} \in \Lambda^* \}  } \{ \sum u_{\delta^*} e^{\delta^*}
\ | \ u_{\delta^*} = 0 \mbox{ for } \ \delta^* \in K  \}.
\ele(CZ*)

By Proposition 2 and Section 7 in \cite{R96} or (\req(sects)) in Appendix\footnote{Indeed, a combinatorial basis for the space of sections of an effective convex divisor class in ${\PZ}_{(\Sigma, L)}$ was given in \cite{R96} Proposition 2, which we will recall in Appendix for convenience here.}, the space of anti-canonical bundle of $\PZ_{\Sigma}$ or $\PZ_{\Sigma_0}$ is given by  
\bea(lllll)
\Gamma ( \PZ_{\Sigma_0} , {\cal O} ( \iota_0^* (-\kappa_0) ) ) &\simeq &\Gamma ( \PZ_{\Sigma} , {\cal O} ( \iota^* (-\kappa) ) ) &  \simeq & \bigoplus \{ {\rm \CZ} {\sf v}^* \ | \ 
{\sf v}^* \in  \triangle^* \cap L^* \ \} \ .
\elea(Ksects)
By (\req(zv*)) in Appendix, the above basis elements ${\sf v}^*$ have the following homogeneous-coordinate representation using $\zeta = (\zeta_\upsilon)_{\upsilon \in {\cal V} ( \triangle )}$ of ${\PZ}_{\Sigma_0}$ in (\req(prin0)) or $z = (z_\delta)_{\delta  \in \Lambda^{(0)}}$ of ${\PZ}_{\Sigma}$ in (\req(prin)): 
\bea(llll)
\zeta^{{\sf v}^*} = \prod_{\upsilon \in {\cal V} ( \triangle )} \zeta_\upsilon^{{\sf v}^*_\upsilon } , &
z^{{\sf v}^*} = \prod_{\delta \in \Lambda^{(0)}} z_{\delta }^{{\sf v}^*_\delta }  & {\rm where} ~ {\sf v}^*_\delta = \langle {\sf v}^*, \delta \rangle + 1 \in \ZZ_{\geq 0} & {\rm for} ~ \delta \in \Lambda^{(0)}.
\elea(zetzv*)
Note that $z^{{\sf v}^*}$ is determined by $\zeta^{{\sf v}^*}$ via the regular map $\widetilde{\varphi}$ in (\req(zzeta)) : $\widetilde{\varphi}^*(\zeta^{{\sf v}^*}) =  z^{{\sf v}^*}$. Furthermore, the above functions define the zero-divisor for the section ${\sf v}^*$ in ${\PZ}_{\Sigma_0}$ and ${\PZ}_{\Sigma}$ respectively,
$$
\sum_{\upsilon \in {\cal V} ( \triangle )} {\sf v}^*_\upsilon e^{\upsilon \dagger} \in D_{\triangle  ~ 0 }^\dagger , ~ ~  \sum_{\delta \in \in \Lambda^{(0)}} {\sf v}^*_\delta e^{\delta \dagger} \in D_{\triangle}^\dagger ,
$$
which are related by $p*$ in (\req(p*)). The zero-loci of a generic section in (\req(Ksects)) define a quasi-smooth hypersurface $X$ of $ \PZ_{\Sigma}$, and the same for $X_0$ in $\PZ_{\Sigma_0}$. Then $X = \varphi^{-1}(X_0)$ with a dominating morphism induced by $\varphi: \PZ_{\Sigma} \longrightarrow \PZ_{\Sigma_0}$:
\be 
\varphi ( = \varphi_{\rm restriction}) : X \longrightarrow X_0 .
\ele(XX0)
Note that both $X$ and $X_0$ are disjoint with 0-dimensional torus-orbits in the toric variety $\PZ_{\Sigma}$ or $\PZ_{\Sigma_0}$ (see Appendix, the case $m_{\tt s} = m = 0$ in Proposition \ref{prop:rhoRf} (III)). By (\req(XX0)), $X$ is a CY $(n-1)$-space  Gorenstein over the (possibly) singular CY space $X_0$. By condition (\req(Lam0)) and the property of reflexive polytope which contains no integral points in the interior except the origin,  $X$ is indeed non-singular when $n \leq 4$, hence a crepant resolution of  $X_0$ (see, e.g., \cite{R89} Proposition 2, or \cite{R96} Proposition 6). We shall write a general element of (\req(Ksects)) in the form
\be
\sum_{w^* \in {\cal V} ( \triangle^* )} w^* + a 0^* + \sum_{{\sf v}^* \in  \triangle^* \cap L^* } \alpha_{{\sf v}^*} {\sf v}^* , ~ ~ a \neq 0, ~ ~ \alpha_{{\sf v}^*} \in \CZ 
\ele(para)
regarded as the hypersurface  deformation of
\bea(ll)
X : & {\tt f} (z) = \sum_{w^* \in {\cal V} ( \triangle^* )} z^{w^*} + a z^{0^*} = 0 ~ \subset \PZ_{\Sigma}  ; \\
X_0 : & f (\zeta) = \sum_{w^* \in {\cal V} ( \triangle^* )} \zeta^{w^*} + a \zeta^{0^*} = 0 ~ \subset \PZ_{\Sigma_0} .
\elea(Margin)
In this section, we are going to determine the deformation space of $X$ through parameters in (\req(para)).  

\subsection{Deformation of anti-canonical hypersurfaces in minimal toric variety\label{ssec.DX0}} 
In this subsection, we determine the deformation space of $X_0$ in the minimal toric variety $\PZ_{\Sigma_0}$ by using homogeneous coordinates $\zeta = (\zeta_\upsilon)_{\upsilon \in {\cal V} ( \triangle )}$. First we will show that by the linear change of variables of $ \CZ_{\widetilde{\Sigma}_0}$ in (\req(prin0)), the parameters $\alpha_{{\sf v}^*}$'s in (\req(para)) can be reduced to those in faces of dimension less than $n-1$ by Jacobian-ring technique, as argued in \cite{R91, R93y} for the simplicial-$\triangle$ case, i.e. Fermat hypersurfaces in weighted projective spaces.
By (\req(zetzv*)), monomials in the expression (\req(Margin)) of $f (\zeta)$  are given by 
$$
\begin{array}{ll}
\zeta^{0^*} = \prod_{\upsilon \in {\cal V} ( \triangle )} \zeta_\upsilon, &
 \zeta^{w^*} = \prod_{\upsilon \in {\cal V} ( \triangle )\setminus F(w^*)} \zeta_\upsilon^{w^*_\upsilon } ~ ~ (w^*_\upsilon \geq 1 ) ~ ~ {\rm for} ~ w^* \in {\cal V} ( \triangle^* ),
\end{array}
$$
where $F(w^*)$ is the facet of $\triangle$ dual to $w^*$; and the coordinates for an element ${\sf x}^* $ in the interior ${\rm Int}(F^*(x))$ of a facet $F^*(x)$ of $\triangle^*$ dual to a vertex $x \in  {\cal V} ( \triangle )$ are expressed by
\be
\zeta^{{\sf x}^*} = \prod_{\upsilon \in {\cal V} ( \triangle ), \upsilon \neq  x } \zeta_\upsilon^{{\sf x}^*_\upsilon } , ~ ~ ~  {\sf x}^*_\upsilon \geq 1 . 
\ele(v*facet)
The Jacobian ring of $f (\zeta)$ in (\req(Margin)) is the ideal in $\CZ [ \zeta ]$ generated by partial derivatives of $f (\zeta)$: 
$J (f) = \langle \partial_\upsilon f (\zeta ) \rangle_{\upsilon \in {\cal V} ( \triangle )}$, where $\partial_\upsilon := \frac{\partial}{\partial \zeta_\upsilon}$. Note that by (\req(Ksects)), $f (\zeta) \in J (f)$ since $ f (\zeta ) = \sum_{ \upsilon \in {\cal V} ( \triangle)} \beta_\upsilon \zeta_\upsilon \partial_\upsilon f (\zeta )$, where $\beta_\upsilon$'s are positive rational numbers satisfying $\sum_{\upsilon \in {\cal V} ( \triangle)} \beta_\upsilon = 1$, and $\sum_{\upsilon \in {\cal V} ( \triangle)} \beta_\upsilon \upsilon = 0$ in $\triangle$.   
\begin{lem}\label{lem:Jac}\footnote{ This lemma and the relation (\req(X0Df)) in this subsection correspond to results in $n$-simplex $\triangle$ case through the Jacobian-ring calculation of a Fermat-type polynomial, which is easily performed (see, e.g. \cite{R91} or \cite{R93y} Theorem 1) , where the coefficients $a_{\ell, \ell'}$'s in the proof of Lemma \ref{lem:Jac} here are zeros when $n \geq 4$, but could be non-zeros in the case $n=3$.} Let $X_0$ be defined by $f (\zeta)$ in (\req(Margin)) with a generic $a$, and ${\sf x}^*$ be a $L^*$-integral element in interior ${\rm Int}(F^*(x))$ of a facet $F^*(x)$ of $\triangle^*$ in (\req(v*facet)). Then 
\bea(ll)
{\sf x}^* \equiv & \sum \{ {\rm \CZ} {\sf v}^* | {\sf v}^* \in L^* \cap (\partial \triangle^* \setminus \bigcup_{F^*:{\rm facet ~of ~ } \triangle^*} {\rm Int}(F^*)) \}   \pmod{J(f)}.
\elea(xJf)
\end{lem}
{\it Proof.} For an element ${\sf x}^*$ in (\req(v*facet)), let $m ( = m ({{\sf x}^*}) ) := {\rm min}\{{\sf x}^*_\upsilon | \upsilon \in {\cal V} ( \triangle ) \setminus \{ x \}  \} \geq 1$, and write $\zeta^{{\sf x}^*}= \eta^{{\sf x}^*} ( \prod_{\upsilon \in {\cal V} ( \triangle )  \setminus  \{ x \} } \zeta_\upsilon )^m $, where $\eta^{{\sf x}^*} = \prod_{\upsilon \in {\cal V} ( \triangle ) \setminus  \{ x \} } \zeta^{{\sf x'}^*} $ with ${\sf x'}^{*}_\upsilon = {\sf x}^*_\upsilon - m \geq 0$. Since $\partial_\upsilon f \in J(f)$, 
$$
\zeta^{{\sf x}^*} = \eta^{{\sf x}^*} (\prod_{\upsilon \in {\cal V} ( \triangle ) \setminus \{ x \}} \zeta_\upsilon)^m \equiv (-a)^m\eta^{{\sf x}^*} (\sum_{w^* \in {\cal V} ( \triangle^* )} \partial_\upsilon \zeta^{w^*})^m = (-a)^m \sum_{ \{ m_{w^*}| w^* \in {\cal V} ( \triangle^* )\}} \eta^{{\sf x}^*}_{\{ m_{w^*}\}} ~ \pmod{J(f)} ,
$$
where the index $\{ m_{w^*}\}$ runs $m_{w^*} \in \ZZ_{\geq 0}$ with $ \sum_{w^* \in {\cal V} ( \triangle^* )}m_{w^*} = m$, and  
$$
\eta^{{\sf x}^*}_{\{ m_{w^*}\}} :=  \eta^{{\sf x}^*} \prod_{w^* \in {\cal V} ( \triangle^* )} (\partial_\upsilon \zeta^{w^*})^{m_{w^*}}  = (\prod_{w^* \in {\cal V} ( \triangle^* )} {w^{*}_\upsilon}^{m_{w^*}})  \prod_{\upsilon \in {\cal V} ( \triangle ) } \zeta_\upsilon^{x^*_{\{ m_{w^*}\}; \upsilon }}.
$$
Note that $\partial_\upsilon \zeta^{w^*} = 0$ when $w^* \in F^*(x)$, (here we use the convention $0^0 := 1$). 
Using $\frac{\partial_\upsilon \zeta^{w^*}}{\partial_\upsilon \zeta^{0^*}} = w^*_\upsilon \frac{\zeta^{w^*}}{\zeta^{0^*}} $ and $\eta^{{\sf x}^*} = \frac{\zeta^{{\sf x}^*}}{(\partial_\upsilon \zeta^{0^*})^m }$, one finds 
$$
\eta^{{\sf x}^*}_{\{ m_{w^*}\}} = \zeta^{{\sf x}^*} \frac{\prod_{w^* \in {\cal V} ( \triangle^* )} ( \partial_\upsilon \zeta^{w^*})^{m_{w^*}}}{(\partial_\upsilon \zeta^{0^*})^m } = (\prod_{w^* \in {\cal V} ( \triangle^* )} {w^{*}_\upsilon}^{m_{w^*}}) \zeta^{{\sf x}^*} \prod_{w^* \in {\cal V} ( \triangle^* )}  \frac{(\zeta^{w^*})^{m_{w^*}}}{(\zeta^{0^*})^{m_{w^*}}}.
$$
Then  $\eta^{{\sf x}^*}_{\{ m_{w^*}\}} \neq 0$ iff $m_{w^*} = 0$ for all $w^* \in F^*(x)$, in which case,  zeros of the monomial $\eta^{{\sf x}^*}_{\{ m_{w^*}\}}$ define a toric divisor of $\PZ_{\Sigma_0}$ linear equivalent to the zero-divisor of $\zeta^{{\sf x}^*}$, equivalently ${\sf x}^*(\{ m_{w^*}\}) (:= {\sf x}^*+ \sum_{w^* \in {\cal V} ( \triangle^* )} m_{w^*} w^* ) \in L^* \cap \triangle^*$ by (\req(Ksects)). By $x^* \in  {\rm Int} ( F^*(x))$  and $0^* \in {\rm Int}(\triangle^*) $, one finds   ${\sf x}^*(\{ m_{w^*}\}) \neq {\sf x}^*$. Therefore 
$$
{\sf x}^* \equiv (-a)^{m({\sf x}^*)} \sum_{\{m_{w^*}| w^* \in {\cal V} ( \triangle^* )\} }  {\sf x}^*(\{ m_{w^*}\}) \prod_{w^* \in {\cal V} ( \triangle^* )} {w^{*}_\upsilon}^{m_{w^*}} ~ \pmod{J(f)}.
$$
where  $m_{w^*} \in \ZZ_{\geq 0}$ and $m_{w^*}=0$ for $w^* \in F^*(x)$   with $ \sum_{w^* \in {\cal V} ( \triangle^* )}m_{w^*} = m({\sf x}^*) (\geq 1)$, ${\sf x}^*(\{ m_{w^*}\}) \in L^* \cap \triangle^* \setminus \{ {\sf x}^* \}$. Denote the vector space on the right side of (\req(xJf)) by $V$, and label elements in $L^* \cap \bigg(\cup_{F^*: {\rm facet ~ of ~} \triangle^*} {\rm Int}( F^*) \bigg)$ by  ${\sf x}_\ell^*$ for $ 1 \leq \ell \leq r$. By $f (\zeta) \in J (f)$, $0^*  \equiv V \pmod{J(f)}$. Hence the above relation implies that for each $\ell$, there exist a positive integer $m({\sf x}_\ell^*)$ and  $a_{\ell, \ell^\prime} \in \ZZ$ for $\ell^\prime \neq \ell$ so that  
$$
{\sf x}_\ell^* - (-a)^{m({\sf x}_\ell^*)} \sum_{\ell^\prime \neq \ell} a_{\ell, \ell^\prime} {\sf x}_{\ell^\prime}^* \equiv V \pmod{J(f)} .
$$
Therefore for generic $a$, ${\sf x}_\ell \equiv V \pmod{J(f)}$ for all $\ell$. This shows the relation (\req(xJf)).
$\Box$ \par  \vspace{.2in} \noindent

By Lemma \ref{lem:Jac}, we may assume $\alpha_{{\sf v}^*}$'s appeared in the parameters of (\req(para)) only for ${\sf v}^*$ in faces of $\triangle^*$ of dimension less than $n-1$, i.e. $v^* = \delta^* \in \Lambda^{*(0)}$ by condition (\req(Lam0)) for the triangulation $\Lambda^*$ of $\partial \triangle^*$. However, there are still constraints (modulus $J(f)$) remained among the coefficients $a$ and $\alpha_{{\sf v}^*}$'s. Since ${\cal V} (\triangle^*) \subseteq \Lambda^{* (0)}$, the monomials $w^*$'s in $f(\zeta)$ can be absorbed into the rest parameters of (\req(para)), expressed in the form of moduli space using variables in (\req(CZ*)):
\be
\zeta^{0^*} + \sum_{\delta^* \in \Lambda^{* (0)}} u_{\delta^*} \zeta^{\delta^*}, ~ u = (u_{\delta^*}) \in \CZ_{\widetilde{\Sigma^*}}.
\ele(moduli)
Note that the above space depends on the triangulation $\Lambda^*$, which might not contain the defining equation of $X_0$ in (\req(Margin)), but with sufficient variables to describe the deformation near $X_0$. We are going to describe a parametrization of (\req(moduli)) using the lattice $D_{\triangle^* }^\dagger$ in (\req(exdual)). Choose a base element $b = (b_{\delta^*})$ in $\CZ_{\widetilde{\Sigma^*} }$. An element $\mu $ in $D_{\triangle^* ~ \QZ}^\dagger$, denoted by
\be
\mu = \sum _{\delta^* \in \Lambda^{* (0)}} \mu_{\delta^*} e^{\delta^* \dagger} \in D_{\triangle^* ~ \QZ}^\dagger,
\ele(mu)
gives rise to an one-parameter family in (\req(moduli)):
\be
X_0 (b ;\mu  ):    \zeta^{0^*}  + s^{\mu} \cdot b = 0, ~ ~  ~ s^{\mu} \cdot b := \sum_{\delta^* \in \Lambda^{* (0)}} s^{\mu_{\delta^*}} b_{\delta^*} \zeta^{\delta^*} ~ ~ (s \in \CZ^*).
\ele(X0sigma)
The families  $X_0 (b ; \alpha \mu )$ and $X_0 (b ; \mu )$ for $\alpha \in \QZ_{>0}$ are related by the change of parameter $s \mapsto s^\alpha$, hence essentially are the same. In particular,  $X_0 (b ; e^{\delta^* \dagger})= 
\zeta^{0^*}  + s_{\delta^*} b_{\delta^*} \zeta^{\delta^*} + \sum_{\delta^{' *} \neq \delta^* \in \Lambda^{* (0)}} b_{\delta^{'*}} \zeta^{\delta^{'*}}$,
by which we introduce the multi-variable deformation,
\be
X_0 (b ; \{ e^{\delta^* \dagger}|e^{\delta^*} \in \Lambda^{(0} \} ):   \zeta^{0^*}  + \sum_{\delta^* \in \Lambda^{*(0)}} s_{\delta^*} b_{\delta^*} \zeta^{\delta^*}= 0 . 
\ele(Xmul)
Note that the above expression provides a parametrization of ${\bf T}(D_{\triangle^*})$-orbit in $\CZ_{\widetilde{\Sigma^*} }$ containing the base point $b$ in  (\req(moduli)), and (\req(X0sigma)) is obtained by the substitution $s_{\delta^*}= s^{\mu_{\delta^*}} $ in (\req(Xmul)). With $b$'s in all ${\bf T}(D_{\triangle^*})$-orbits, (\req(Xmul)) provides a simultaneous parametrization of the deformation space (\req(moduli)). However among families in (\req(X0sigma)), 
those induced from automorphisms of toric variety $\PZ_{\Sigma_0}$ are redundant ones. Indeed by (Remark of) Lemma \ref{lem:aut} in Appendix, an one-parameter group of ${\bf T}(L)$, i.e.  $\beta^{\prime *}(\xi) $ in (\req(exdual)) for $\xi \in L$,  acts on ${\sf v}^*$ in (\req(Ksects)) by $(s,  {\sf v}^*) \mapsto s^{\langle {\sf v}^* , \xi \rangle} {\sf v}^*$, where $s \in \CZ^* (=\CZ/\ZZ)$ is the independent variable of one-parameter group $\xi$ in the expression of $X_0 (b ;\beta^{\prime *} (\xi) )$ (\req(X0sigma)) with $\beta^{\prime *}(\xi)= \sum _{\delta^* \in \Lambda^{* (0)}} \langle \delta^*, \xi \rangle e^{\delta^* \dagger}$. Any two deformations in (\req(moduli)) differ by $X_0 (b ;\beta^{\prime *}(\xi) )$ for $\xi \in L$ are equivalent up to a family of transformations of $\PZ_{\Sigma_0}$-automorphism. Hence deformations of $X_0$ modulus $\PZ_{\Sigma_0}$-automorphism are identified with  $D_{\triangle^* ~ \QZ}^\dagger$ modulus $L_{\QZ}$, which is the same as  ${\bf n}_{\triangle^* ~ \QZ}^\dagger$ by 
(\req(exdual)). Hence we obtain the following description of deformation classes: 
\bea(lll)
{\rm Def}(X_0)_{\CZ} & \simeq & {\bf n}^\dagger_{\triangle^*~\CZ}.
\elea(X0Df)
Note that ${\bf n}^\dagger_{\triangle^*~\CZ}$ is defined through integral points in  the faces of $\triangle^*$, which does not depend on the detailed triangulation structure of $\Lambda^*$. A realization of ${\rm Def}(X_0)_{\CZ}$ is provided by a lifting of ${\bf n}^\dagger_{\triangle^*~\CZ}$ in $D_{\triangle^* ~ \CZ}^\dagger$ through $\iota^{\prime *}$ in (\req(exdual)). 

\subsection{Deformation of anti-canonical hypersurface in toric variety defined by reflexive polytope\label{ssec.DX}} 
In this subsection, we describe a combinatorial representation of the deformation space of  anti-canonical hypersurface $X$ of $\PZ_{\Sigma}$ in (\req(Margin)). By (\req(K3)) and (\req(PiDe)), ${\rm Def}(X)_{\CZ}$ can be regarded a subspace of the Hodge space $H^{n-2,1}(X)$ in $H^{n-1}(X)$, similarly  ${\rm Def}(X_0)_{\CZ} \subseteq H^{n-2,1}(X_0) \subseteq H^{n-1}(X_0)$. Using Mayer-Vietoris cohomology sequences, one finds the cohomology relation between $X$ and $X_0$ through $\varphi$ in (\req(XX0)):
$$
\begin{array}{ll}
H^{n-1}(X) \simeq H^{n-1}(X_0) \oplus \bigg(\oplus_E H^{n-1}(E)\bigg) , &  H^{n-2, 1}(X) \simeq H^{n-2, 1}(X_0) \oplus \bigg(\oplus_E H^{n-3, 0}(E)\bigg)
\end{array}
$$
where $E$'s run the exceptional divisors in $X$, and $H^{n-1}(X)$ is related to $H^{n-1}(E)$ by restriction of cohomology class, equivalent to the cup-product of $(1,1)$-form representing the Chern class of ${\cal O}(E)$, (see, e.g. \cite{R90}  Theorem 2). By Proposition \ref{prop:rhoRf} (III) (with $\rho = - \kappa, m_{\tt s} = n-1$) in Appendix, $H^{n-3,0}(E)=0$ except $\varphi (E)$ is a hypersurface of a 
$(n-2)$-dimensional toric subvariety of $\PZ_{\Sigma_0}$, in which case $E = X \cap \overline{\rm O}_{\tt s} $ for a vertex ${\tt s} \in \Lambda^{(0)} \cap {\rm Int}(F)$ with $F$ a 1-dimensional face of $\partial \triangle$, and $H^{n-3, 0}(E)$ is represented by (\req(KXb)) with $F^*_\rho = F^*, m=n-2$. Indeed, the exceptional divisor $E$ is generic $\PZ^1$-bundle over $\varphi (E)$, and the contribution $H^{n-3, 0}(E)$ in the deformation space are induced from canonical sections of $\varphi (E)$. Hence ${\rm Def}(X)_{\CZ} $ is the direct sum of ${\rm Def}(X)_{0 ~\CZ}$ with all those $H^{n-3, 0}(E)$'s.  By (\req(X0Df)), we obtain the combinational representation of deformation space of $X$:
\be
{\rm Def}(X)_{\CZ} \simeq {\bf n}^\dagger_{\triangle^*~\CZ} \oplus \bigoplus_{F^*, ~ \nu_{F^*} } \CZ ~ \nu_{F^*} 
\ele(DefX)
where $F^*$ runs over codim. 2 faces of $\triangle^*$, and  $\nu_{F^*} \in ( {\rm Int}(F^*) \cap L^* ) \times ({\rm Int}(F) \cap L)$. For $\mu \in D_{\triangle^* ~ \QZ}^\dagger$ in (\req(mu)), the deformation $\iota^{\prime *} \mu \in {\bf n}^\dagger_{\triangle^*~\CZ}$ in ${\rm Def}(X)_{\CZ}$ is induced from $X_0 (b ;\mu  )$ in (\req(X0Df)), hence by (\req(Ksects)), is represented by 
\be
X (b ;\mu  ):    z^{0^*}  + s^{\mu} \cdot b = 0, ~ ~  ~ s^{\mu} \cdot b := \sum_{\delta^* \in \Lambda^{* (0)}} s^{\mu_{\delta^*}} b_{\delta^*} z^{\delta^*} ~ ~ (s \in \CZ^*).
\ele(Xsigma)
Note that ${\rm Def}(X)_{\CZ}$ in (\req(DefX)) depends only on the 0-skeleton $\Lambda^{* (0)}$ of $\Lambda^*$. The structure in the deformation space relevant to the detailed triangulation $\Lambda^*$ is the degeneration cone of $X$ discussed in the next subsection, similar to the relation between (\req(PicX)) and Kahler cone in Subsection \ref{ssec.Kcone}.

\subsection{ Degeneration cone in moduli space of anti-canonical hypersurface of a toric variety \label{ssec.+deg}}
In this subsection, we discuss the degeneration of anti-canonical hypersurfaces of a toric variety in the moduli space (\req(moduli)), whose structure  depends on the simplicial decomposition $\Lambda^*$ of $\partial \triangle^*$. 
In (\req(Xmul)), we use the one-parameter families of $e^{\delta^*}$'s to parameterize the moduli space (\req(moduli)). However, the structure of ${\bf n}^\dagger_{\triangle^*~\QZ}$ in (\req(DefX)) will be better represented by a $D^\dagger_{\triangle^* ~ \QZ}$-basis  dual to a basis of $D_{\triangle^* ~ \QZ}$ which contains a ${\bf n}_{\triangle^*~\QZ}$-basis.  For $(n-1)$-simplex ${\tt s}^* \in \Lambda^{* (n-1)}$, let $\{ e^{\delta^*} \}_{\delta^* \in {\tt s}^* \cap \Lambda^{* (0)}} \cup \{ n^{\delta^{* \prime }}\}_{\delta^{* \prime} \in  \Lambda^{* (0)} \setminus {\tt s}^*}$ be the basis of $D_{\triangle^* ~ \QZ}$ in Lemma \ref{lem:locS} with respective to $(L^*, \Lambda^*)$, and  $\{ e_{\delta^* , {\tt s}^*} \}_{\delta^* \in {\tt s}^* \cap \Lambda^{* (0)}} \cup \{ n_{\delta'^*, {\tt s}^*}\}_{\delta^{* \prime} \in  \Lambda^{* (0)} \setminus {\tt s}^*}$ is the dual basis. As in (\req(rhos)), any element $\mu \in D_{\triangle^* ~ \QZ}^\dagger$ can be expressed  by 
\be
\mu =  \mu_{{\tt s}^* }^\prime + \mu_{{\tt s}^*}, ~ ~\mu_{{\tt s}^* }^\prime = \sum_{\delta^* \in {\tt s}^*} \langle \mu  , e^{\delta^*} \rangle e_{\delta^* , {\tt s}^*} \in \beta^{\prime *}(L_{\QZ}) , ~ ~ \mu _{{\tt s}^*} = \sum_{\delta^{* \prime} \notin {\tt s}^* }\langle \mu  , n^{\delta^{* \prime}} \rangle n_{\delta^{* \prime}, {\tt s}^*} \in D^\dagger_{{\tt s}^*}.
\ele(mus)
Parallel to (\req(Pcone)), Proposition \ref{prop:+cone} and (\req(n+)), we define the positive cone in $D_{\triangle^* ~ \QZ}^\dagger$, 
$$
C^\dagger_{* ~ +} ~ (= C^\dagger_{\Lambda^*, +}) = \{ \mu \in  D_{\triangle^* ~ \QZ}^\dagger | \langle n^{\delta^{* \prime}}_{{\tt s}^*}, e^{\delta^{* \prime}} \rangle  ~ (= \langle \mu_{\tt s}, e^{\delta^{* \prime}} \rangle ) ~ \geq 0 ~ ~ (\delta^{* \prime} \notin {\tt s}^*  \in \Lambda^{* (n-1)})  \},
$$
where $n^{\delta^{* \prime}}_{{\tt s}^*} \in {\bf n}_{\triangle^* ~ \QZ}$ are defined by (\req(nsdel)) for a 
$(n-1)$-simplex ${\tt s}^* \subset L^*_{\RZ}$.  All $n^{\delta^{* \prime}}_{{\tt s}^*}$'s generate the cone ${\bf n}_{\Lambda^*, +}=\sum_{\delta^{* \prime} \notin{\tt s}^* \in \Lambda^{* (n-1)}} \QZ_{\geq 0} ~ n^{\delta^{* \prime}}_{{\tt s}^*}$ in ${\bf n}_{\triangle^* ~ \QZ}$.   As in (\req(ndual)), we define the degeneration cone
\bea(ll)
{\cal C}_{\rm Deg} (X) (= {\cal C}_{{\rm Deg}, \Lambda^*}(X) )= {\bf n}^\dagger_{\Lambda^*, +} := {\rm dual ~ cone ~ of ~ } {\bf n}_{\Lambda^*, +}  &\subseteq  {\bf n}^\dagger_{\triangle^*~ {\rm \QZ}},
\elea(degc)
which consists of non-negative (one-parameter) degeneration classes. The interior of $C^\dagger_{* ~ +} $ consists of all positive degeneration of anti-canonical hypersurfaces in $\PZ_{\Sigma}$. In particular, elements in ${\rm Int}(C^\dagger_{* ~ +}) \cap {\rm Int}(D_{\triangle^*~ +}^{\dagger})$ are CY degenerations with $(s=0)$-limit $z^{0^*}$, called the maximal unipotent degenerations. Note that any such deformation $\mu$ is equivalent to $\mu _{{\tt s}^*}$ in (\req(mus)) for ${\tt s}^* \in \Lambda^{(n-1)}$, which is a degeneration with $(s=0)$-limit $\sum_{{\sf v}^* \in \Lambda^{*(0)} \cap {\tt s}^*} z^{{\sf v}^*}$.

\section{Mirror Symmetry of Anti-canonical Hypersurfaces in Toric Varieties of Reflexive Polytopes \label{sec.mir}}
\setcounter{equation}{0}
We start with a triangulation $\Lambda$ of a reflexive polytope $(\triangle, L)$ and a triangulation $\Lambda^*$ of the dual polytope $(\triangle^*, L^*)$. In Sections \ref{sec.PG} and \ref{sec.Defm}, we derive the combinatorial representation of ${\rm Pic}(X)_{\CZ}$ and ${\rm Def}(X)_{\CZ}$ for a generic anti-canonical hypersurface $X$ of $\PZ_{\Sigma(\Lambda)}$ in (\req(PicX)) and (\req(DefX)) respectively, then find an explicit form of Kahler cone  ${\bf n}^\dagger_{\Lambda, +}$ of $\PZ_{\Sigma(\Lambda)}$ in (\req(Kcone)), and degeneration cone ${\bf n}^\dagger_{\Lambda^*, +}$ in the moduli space (\req(moduli)) of $X$ in (\req(degc)). On the other hand, $(\triangle^*, L^*)$ gives rise to a ${\bf T}(L^*)$-toric variety $\PZ_{\Sigma(\Lambda^*)}$, whose anti-canonical global sections possess a basis represented  by $\triangle \cap L$ as in (\req(Ksects)). Let $X^*$ be a generic anti-canonical hypersurface of $\PZ_{\Sigma(\Lambda^*)}$. By regarding $L^*$ as the group of one-parameter subgroups of ${\bf T}(L^*)$ with the character group $L$, one finds the combinatorial representatives of ${\rm Pic}(X^*)_{\CZ}$ are the same as ${\rm Def}(X)_{\CZ}$ in (\req(DefX)), and a similar  identification also for ${\rm Def}(X^*)_{\CZ}$ and ${\rm Pic}(X)_{\CZ}$.  Therefore $X$ and $X^*$ constitute a mirror CY pair so that the relation (\req(PiDe)) holds. Furthermore, the Kahler cone ${\cal C}_{{\rm Pic}, \Lambda^*}(X^*)$ of $\PZ_{\Sigma(\Lambda^*)}$ is identified with the degeneration cone ${\cal C}_{{\rm Deg}, \Lambda^*} (X)$ for $X$ in (\req(degc)), and the same for the equality between ${\cal C}_{{\rm Deg}, \Lambda} (X^*)$  and  the Kahler cone ${\cal C}_{{\rm Pic}, \Lambda} (X) $ of $\PZ_{\Sigma(\Lambda)}$. In summary, the following structures between a mirror anti-canonical hypersurfaces of two toric varieties defined by a reflexive and its dual polytopes are interchangeable by identifying their toric representatives:
$$
\begin{array}{lll}
{\rm Picard ~ group :~ Pic_{\CZ}} & \Longleftrightarrow & {\rm Deformation ~ space: ~ Def_{\CZ}} , \\
{\rm A ~line ~bundle} & \Longleftrightarrow & {\rm One-parameter~ deformation} , \\
{\rm Principal~bundle~space ~ over ~ toric ~ variety} & \Longleftrightarrow & {\rm Moduli~space~of~ anti-canonical ~ hypersurace} , \\
{\rm Kahler ~ cone}: ~ {\cal C}_{\rm Pic} & \Longleftrightarrow & {\rm Degeneration ~ cone}:~ {\cal C}_{\rm Deg},  \\
{\rm A ~positive~effective~ divisor~ class} & \Longleftrightarrow & {\rm a ~maximal ~unipotent~degeneration} . \\
\end{array}
$$
The ${\cal C}_{{\rm Pic}, \Lambda}(X)$ of $X$ in above describes the Kahler cone of $\PZ_{\Sigma(\Lambda)}$, which is equal to degeneration cone ${\cal C}_{{\rm Deg}, \Lambda}(X^*)$ of its mirror $X^*$, both depending on the triangulation $\Lambda$ of $\partial \triangle$. In particular for a CY 3-fold $X$, there is the flop process for certain $\PZ_{\Sigma(\Lambda)}$ as described in Example 3 of Section \ref{ssec.Kcone}, where $X$ is birational, but not biregular, to another CY 3-fold $X^\natural$ with the same combinatorial representation of Picard group (\req(PicX)) and deformation space (\req(DefX)).   $X$ and $X^\natural$ also share the same degeneration cone (\req(degc)), but with different Kahler cones by ${\rm Int}({\cal C}_{{\rm Pic}}(X)) \cap {\rm Int}({\cal C}_{{\rm Pic}}(X^\natural)) = \emptyset$. The mirror $X^*$ and $X^{\natural ~ *}$ of $X$ and $X^\natural$  respectively are biregular, but their degeneration cones are different with ${\rm Int}({\cal C}_{{\rm Deg}}(X^*)) \cap {\rm Int}({\cal C}_{{\rm Deg}}(X^{\natural ~ *})) = \emptyset$, which serves the mirror of flop of CY 3-folds.

\section{Concluding Remarks \label{sec.F}}
From the framework of toric geometry, we find a combinatorial  representation of Picard group and deformation space of anti-canonical  hypersurfaces $X$ of a toric variety defined by reflexive polytope. By identifying the toric representatives of the representation, we establish the cohomology correspondence between the mirror pair previously found in \cite{B}, as an extension of the mirror correspondence in \cite{R91}. This paper can be regarded as a continuation of the work in \cite{R96} about Picard group of hypersurfaces in toric varieties. Here we have further identified the structure of Kahler cone and degeneration cone of mirror CY hypersurfaces using toric techniques. As a consequence, a different parametrization of the moduli space in CY degeneration provides the mirror notion equivalent to flops in CY 3-folds. For the purpose of exploring the essential role of simplicial-cone structure of toric variety in mirror symmetry, we concentrate here only on the positive line bundles of $X$ induced from the Kahler cone of the ambient toric variety. However in the most general situation, the complete Kahler cone or positive deformation cone of $X$  should add the extra contribution of $\nu_F$'s in (\req(PicX)) or $\nu_{F^*}$'s in (\req(DefX)), into (\req(Kcone)) or (\req(degc)), the relation of which not immediately apparent, and not yet studied in the present article. These structures are essential for the understanding of quantum cohomology product in mirror symmetry. Work along these lines is under consideration. Here, just to keep things simple, we restrict our attention only on the mirror cone-structure with no contribution from $\nu_F$'s or $\nu_{F^*}$'s. We leave the further discussion of our results, and possible generalizations to future work.

\section*{Appendix: Some Basic Facts in Toric Geometry}
\setcounter{equation}{0}
In this Appendix, we summarize results in \cite{R96} and some basic facts in toric geometry used in this paper for easy reference.  Some should be well-known (but we could not find an explicit reference). For completeness, we include some technical arguments for certain statements here. Results in this section hold for  a general polytope $(\triangle, L)$ with a simplicial decomposition $\Lambda$ of $\partial \triangle$ with $\Lambda^{(0)} \subset L \cap \partial \triangle$, i.e. condition (\req(Lam0)) or (\req(L0)) not necessary required. For convenience of notions, we consider the cases when  the rational toric divisor class ${\cal O}(\iota^* \rho)$ over ${\PZ}_{\Sigma}$ is effective and convex, i.e.  
\be
\rho = \sum_{\delta \in \Lambda^{(0)}} \rho^\delta e^{\delta \dagger} \in D_{\triangle}^\dagger  ~ ~ (  \rho^\delta \in \ZZ_{>0} ), \tag{A1}
\ele(rho)
and the piecewise linear functional on $L_{\RZ}$, linear on every simplicial cone $\sigma ({\tt s})$ of $\Sigma$ with value $\rho^\delta$ on $\delta \in \Lambda^{(0)}$, satisfies the convex property \cite{O88, OP, R96}, a condition equivalent to the non-negative condition in Section \ref{ssec.Kcone} of this paper. Define the $\rho$-dual polytope of $\triangle$ by\footnote{The $\triangle^*_\rho$ here is denoted by $\triangle (L^*)_\rho$ in \cite{R96} (24). }   
\be
\triangle^*_\rho = \{ y \in  L^*_{\RZ} \ | \ \langle y, \delta \rangle \ \geq -\rho^\delta \  \ \mbox{for } \ \delta \in \Lambda^{(0)} \}. \tag{A2}
\ele(t*rho)
Note that the effective-convex condition of $\rho$ is equivalent to the convex-cone property of $\triangle^*_\rho$ (see, e.g. Proposition 1 and (24) in \cite{R96}). By the convex property of $\rho$, each $m$-face $F_\rho^*$ of $\partial \triangle^*_\rho$ for $0 \leq m \leq n-1$ gives rise to a $(n-1-m)$-polytope $F_\rho$ in $\partial \triangle$ generated by finite elements in $L_{\QZ}$. We shall assume $F_\rho$ is generated by a finite number of elements in $\Lambda^{(0)} \cap F$ for some $(n-1-m)$-face $F$ of $\partial \triangle$, so that all $F_\rho$'s form a polytope refinement of the face-decomposition of $\partial \triangle$. For convenience, $F_\rho$ will be called a $\rho$-face of $\partial \triangle$, and the set of 0-dimensional $\rho$-faces will be denoted by ${\cal V}_\rho ( \triangle )$.  In general, the $\rho$-face $F_\rho$ may not be equal to the face $F$ (for the details, see the relation between $\rho$-graph-cone ${\cal C}_\rho$ and its dual cone ${\cal C}^*_\rho$ in \cite{R96}  Section 3). Similar to the construction of the minimal toric variety $\PZ_{\Sigma_0}$, there exists the "$\rho$-minimal" toric variety $\PZ_{\Sigma_{\rho ~ 0}}$ for a convex divisor $\rho$ in (\req(rho)), where $\Sigma_{\rho ~ 0}$ is the complete polyhedral cone of $L_{\RZ}$ induced by the $\rho$-face-decomposition of $\partial \triangle$. Since $\Lambda$ is a triangulation of $\rho$-faces and $\rho$ is determined by those $\rho^\delta$'s for $\delta \in {\cal V}_\rho ( \triangle )$ in (\req(rho)), the toric variety $\PZ_{\Sigma}$ dominates $\PZ_{\Sigma_{\rho ~ 0}}$:
\bea(ll)
\varphi (= \varphi_\rho) : \PZ_{\Sigma (\Lambda )} \longrightarrow \PZ_{\Sigma_{\rho ~ 0}} , & {\cal O} (\iota^* \rho) = \varphi^* {\cal O} (\iota_0^* \rho_0) , \tag{A3}
\elea(domi)
where $\rho_0$  is the toric divisor of $\Sigma_{\rho ~ 0}$ defined by $\rho_0:= \sum_{v \in {\cal V}_\rho ( \triangle )} \rho^v e^{v \dagger}$. Hence sections of ${\cal O} (\iota^* \rho)$ over $\PZ_{\Sigma (\Lambda )}$ are the $\varphi$-pull-back of ${\cal O} (\iota_0^* \rho_0)$-sections over $\PZ_{\Sigma_{\rho ~ 0}}$. Note that when $\rho= -\kappa$ in (\req(kappa)), $\triangle^*_{-\kappa}$ is the same as the dual polytope $\triangle^*$ of $\triangle$, and $(- \kappa)$-face $F_{-\kappa}$ is equal to the face $F$ with  
$\PZ_{\Sigma_{\rho ~ 0}} = \PZ_{\Sigma_0}$,  $\rho_0 = - \kappa_0$ in (\req(kappa0)).

Now we describe a combinatorial basis of sections of ${\cal O} (\iota^* \rho)$. First we
consider  the case when all elements in $\Lambda^{(0)}$ are primitive in $L$.  By Proposition 2 in \cite{R96}\footnote{In \cite{R96}, we assume the polytope $\triangle$  contains the origin as the only lattice point in its interior. Indeed by the same argument, the results in Proposition 2 and Theorem 2 of \cite{R96} are still valid for integral polytopes discussed in the appendix here.}, a basis of global sections of ${\cal O}(\iota^* \rho)$ over ${\PZ}_{\Sigma}$ consists of a finite number of integral elements in $\triangle^*_\rho $ as follows:
\bea(llll)
\Gamma ( \PZ_{\Sigma} , {\cal O} ( \iota^* \rho ) ) & \big( \simeq \Gamma ( \PZ_{\Sigma_{\rho ~ 0}} , {\cal O} ( \iota_0^* \rho_0 ) ) \big) &  \simeq & \bigoplus \{ {\rm \CZ} {\sf v}^* \ | \ 
{\sf v}^* \in  \triangle^*_{\rho} \cap L^* \ \} \ .  \tag{A4}
\elea(sects)
Indeed by the exact sequences (21) and (22) in \cite{R96}, the above basis elements ${\sf v}^*$'s are identified with  monomial functions of $\CZ_{\widetilde{\Sigma}}$ in (\req(CZ)) expressed by
the homogeneous coordinates $z = (z_\delta)_{\delta }$ of ${\PZ}_{\Sigma}$ as
\be
z^{{\sf v}^*} = \prod_{\delta \in \Lambda^{(0)}} z_{\delta }^{{\sf v}^*_\delta } ,  ~ ~ ~ {\sf v}^*_\delta:= \langle {\sf v}^*, \delta \rangle + \rho^\delta \in \ZZ_{\geq 0}, \tag{A5}
\ele(zv*)
which are characterized as monomials $z^k ~  (k = \sum_{\delta} k^\delta e^{\delta \dagger} \in D^\dagger_\triangle, k^\delta \geq 0)$ with the "degree" property: $\langle k , \iota (n) \rangle = \langle \rho, \iota (n) \rangle$ for $n \in {\bf n}_\triangle $. In particular, the zero-section $0^*$ corresponds to $z^{0^*} = \prod_{\delta \in \Lambda^{(0)}} z_{\delta }^{\rho^\delta }$ with $\rho$ as its zero-divisor. Furthermore, as rational functions of the fiber space in (\req(prin)), $\frac{z^{{\sf v}^*}}{z^{0^*}}$ is identified with the character function of ${\bf T}(L)$ induced from ${\sf v}^*$, $\chi_{{\sf v}^*} : {\bf T}(L) \longrightarrow \CZ^* (=\CZ/\ZZ)$, regarded as a rational function of  $\PZ_{\Sigma}$.
Therefore through the expression (\req(zv*)), the space of sections in (\req(sects)) is identified with a subspace of polynomial functions of $\CZ_{\widetilde{\Sigma}}$. It is known that the translation $m_{\overline{u}}$ of ${\bf T}(L)$ by an element $\overline{u} \in {\bf T}(L)$ extends an automorphism of $\PZ_{\Sigma}$  which preserves all toric divisors, hence the line bundle ${\cal O}(\iota^*\rho)$. Furthermore, one may lift the translation $m_{\overline{u}}$ to the bundle ${\cal O}(\iota^*\rho)$ fixing the section $z^{0^*}$. Indeed, we have the following result:
\begin{lem}\label{lem:aut}
The ${\bf T}(L)$-action of ${\rm \PZ}_{\Sigma}$ can be lifted to the line bundle ${\cal O}(\iota^*\rho)$ so that ${\bf T}(L)$ acts on the space (\req(sects)) by 
$$
(\overline{u} , {\sf v}^*) \mapsto \chi_{{\sf v}^*}(\overline{u}) {\sf v}^*  , ~ ~ \overline{u} \in {\bf T}(L) , {\sf v}^* \in  \triangle^*_{\rho} \cap L^* ,  
$$
where $\chi_{{\sf v}^*}$ is the character function of ${\bf T}(L)$ determined by ${\sf v}^*$.
\end{lem} 
{\it Proof.} By the definition of ${\cal O}(\iota^*\rho)$ and (\req(zv*)), it suffices to construct a lifting of the ${\bf T}(L)$-action of ${\rm \PZ}_{\Sigma}$ to the fiber space $\CZ_{\widetilde{\Sigma}}$ of (\req(prin)) which fixes the monomial $z^{0^*}$. Let $\delta^j ~ (1 \leq j \leq n)$ be $n$ elements in $\Lambda^{(0)}$ which generate a $(n-1)$-simplex contained in a $\rho$-facet $F_\rho$ of $\triangle$.  Then $\delta^j$'s form a basis of $L_{\QZ}$ with the dual basis $\delta_j \in L^*_{\QZ}$, $\langle \delta_j, \delta^k \rangle = \varepsilon_j^k$.  For $\delta' \in \Lambda^{(0)} \setminus \{ \delta^j \}_{j=1}^n$, one has 
\bea(ll)
\delta' = - \sum_{j=1}^n a^{\delta'}_j \delta^j   & {\rm for}~  a^{\delta'}_j \in \QZ, \tag{A6}
\elea(d'r)
by which we define $n^{\delta'} = e^{\delta'} + \sum_{j=1}^n a^{\delta'}_j e^{\delta^j} \in {\bf n}_{\triangle ~ \QZ}$. Then $n^{\delta'}$'s form a basis of ${\bf n}_{\triangle ~ \QZ}$, and together with $
\{ e^{\delta^j} \}_{j=1}^n$, they form a basis of $D_{\triangle ~ \QZ}$. The basis $
\{ e_{\delta^j} \}_{j=1}^n \cup \{ n_{\delta'} \}_{\delta' \in \Lambda^{(0)}, \delta' \neq \delta^j }$ of $D^\dagger_{\triangle ~ \QZ}$ dual to $\{ e^{\delta^j} \} \cup \{ n^{\delta'} \}$  is related to $e^{\delta \dagger}$'s in (\req(D*)) by 
\bea(lll)
e_{\delta^j } = e^{\delta^j \dagger} - \sum_{\delta' \in \Lambda^{(0)} \setminus \{ \delta^j \}} a_j^{\delta'} e^{\delta'' \dagger}, & 
n_{\delta'} = e^{\delta' \dagger}  ,  &
e^{\delta^j \dagger} = e_{\delta^j }  + \sum_{\delta' \in \Lambda^{(0)} \setminus \{ \delta^j \}} a_j^{\delta'} n_{\delta'} . \tag{A7}
\elea(dnrho)
Under $\beta^*$ in (\req(exa*)), one has $\beta^*(\delta_j) = e_{\delta^j }$ for $1 \leq j \leq n$. We now use the basis $\delta_j$'s of $L^*_{\QZ}$ to express an element ${\sf w}^* \in L^*$  by ${\sf w}^* = \sum_{j=1}^n w^j \delta_j  $ for $w^j \in \QZ$, and write $\overline{u} = u + L \in {\bf T}(L) (=L_{\CZ}/L)$ with $u = \sum_{j=1}^n u_j \delta^j \in L_{\CZ}$ for $u_j \in \CZ$. Then the character $\chi_{\sf w}^*$ of ${\bf T}(L)$ is expressed by $\chi_{{\sf w}^*}(\overline{u}) = \prod_{j=1}^n e^{2 \pi {\rm i}  u_j w^j}$. In particular, $\chi_{\delta_j}(\overline{u}) = e^{2 \pi {\rm i}  u_j}$,  hence $\chi_{\delta_j}(\overline{u}) e_{\delta_j} = e^{2 \pi {\rm i}  u_j} e_{\delta_j}$ for $1 \leq j \leq n$. Since $z^{0^*}$ is identified with $\rho \in D^\dagger_{\triangle}$, expressed by $
\rho = \sum_{j=1}^n \rho^{\delta^j} e_{\delta^j}+ \sum_{\delta' \in \Lambda^{(0)} \setminus \{ \delta^j \}} (\rho^{\delta'} + \sum_{j=1}^n \rho^{\delta^j}  a_j^{\delta'}) n_{\delta'}$ in terms of the basis $\{ e_{\delta^j} \} \cup \{ n_{\delta'} \}$ , 
a ${\bf T}(L)$-lifting of ${\rm \PZ}_{\Sigma}$ to $\CZ_{\widetilde{\Sigma}}$ in (\req(prin)) fixing $z^{0^*}$ is equivalent to a collection of $(d-n)$ linear functionals, $v_{\delta'} \in {\rm Hom}( L_{\CZ},\CZ)$ for $\delta' \in \Lambda^{(0)} \setminus \{ \delta^j \}$, satisfying the equation
$$
\begin{array}{ll}
\sum_{j=1}^n \rho^{\delta^j} u_j + \sum_{\delta' \in \Lambda^{(0)} \setminus \{ \delta^j \}} (\rho^{\delta'} + \sum_{j=1}^n \rho^{\delta^j}  a_j^{\delta'}) v_{\delta'}(u) & = 0,   
\end{array}
$$ 
which by (\req(dnrho)), the ${\bf T}(L)$-action of $\CZ_{\widetilde{\Sigma}}$, $(\overline{u}, (z_{\delta})_{\delta \in \Lambda^{(0)}}) \mapsto (z_{\delta}(\overline{u}))_{\delta \in \Lambda^{(0)}}$, is described by
\bea(ll)
z_{\delta^j}(\overline{u}) = e^{2 \pi {\rm i}(\rho^{\delta^j} u_j + \sum_{\delta'' \notin \{ \delta^k \} } a_j^{\delta''} v_{\delta''}(u))} z_{\delta^j}, ~ ~ z_{\delta'}(\overline{u}) = e^{2 \pi {\rm i}  v_{\delta'}(u)} z_{\delta'} ~ (\delta' \in \Lambda^{(0)} \setminus \{ \delta^j \}).  \tag{A8}
\elea(LCa)
By (\req(d'r)), the vanishing condition, $\rho^{\delta'} + \sum_{j=1}^n \rho^{\delta^j}  a_j^{\delta'}=0$, is equivalent to $\delta'$ lies in the hypersurface of  $L_{\RZ}$ containing the facet $F_\rho$. For a given element $\delta^0 \in \Lambda^{(0)}$ outside the hypersurface containing $F_\rho$, an action in (\req(LCa)) is given by the following $v^{\delta'}$'s:
$$
\begin{array}{ll}
v_{\delta^0}(u) = - \sum_{j=1}^n \frac{\rho^{\delta^j}}{\rho^{\delta^0} + \sum_{j=1}^n \rho^{\delta^j}  a_j^{\delta^0}} u_j  , & v_{\delta'}(u) = 0 ~ {\rm for} ~ \delta' \neq \delta^0 ,
\end{array}
$$ 
which provides a ${\bf T}(L)$-lifting to the line bundle ${\cal O}(\iota^*\rho)$. 
$\Box$ \par \noindent 
{\bf Remark.} As in (\req(prin)) and (\req(zv*)), there is a homogeneous coordinate system $\zeta =(\zeta_\upsilon)_{\upsilon \in {\cal V}_\rho ( \triangle )}$ for $\PZ_{\Sigma_{\rho ~ 0}}$, by which the section ${\sf v}^*$ of ${\cal O}(\iota_0^*\rho_0)$ in (\req(sects)) is identified with $\zeta^{{\sf v}^*}$ (for the case $\rho= -\kappa$, see (\req(prin0)) and (\req(zetzv*))). By the same argument for the expression $\zeta^{{\sf v}^*}$, one can show that Lemma \ref{lem:aut} also holds for the line bundle ${\cal O}(\iota_0^*\rho_0)$ over ${\rm \PZ}_{\Sigma ~ 0}$.
$\Box$ \par \vspace{.2in} \noindent
Suppose $X$ is a quasi-smooth hypersurface $X$ of $\PZ_{\Sigma}$ defined by zeros of a section in (\req(sects)). By using the standard residual argument in algebraic geometry, a basis for sections of the canonical bundle ${\cal K}_X$ of $X$ can be derived from those in (\req(sects)) as follows:
\bea(lll)
H^{n-1, 0} (X) (= \Gamma ( X , {\cal K}_X )  ) & \simeq  & \bigoplus \{{\rm  \CZ} {\sf w}^* \ | \ {\sf w}^* \in  {\rm Int} (\triangle^*_{\rho}) \cap L^* \ \} .  \tag{A9}
\elea(HX)

We now consider the general case for an integral polytope $(\triangle, L)$ with a simplicial decomposition $\Lambda$ of $\partial \triangle$ with $\Lambda^{(0)} \subset L \cap \partial \triangle$. For each $\rho \in \Lambda^{(0)}$, let  $m_\rho$ be the positive integer  so that $\overline{\rho} ~ (:=\frac{\rho}{m_\rho})$ is a primitive element in $L$. Instead of $\rho$'s, the primitive toric divisors in $\PZ_{\Sigma} $ are represented by $\overline{\rho}$'s. An integral toric divisor is an element in the lattice $\overline{D}^\dagger$ dual to $\overline{D}$, where
$$
\begin{array}{lll}
\overline{D} := \bigoplus_{\delta \in \Lambda^{(0)}} \ZZ e^{\overline{\rho}} \subset D_{\triangle ~\QZ} , & ~ & (e^{\overline{\rho}} = \frac{ e^{\rho}}{m_\rho}) , \\
\overline{D}^\dagger := \bigoplus_{\delta \in \Lambda^{(0)}} \ZZ e^{\overline{\rho} \dagger} \subset D^\dagger_{\triangle ~\QZ}, & ~ & (e^{\overline{\rho} \dagger} = m_\rho e^{\rho \dagger}),
\end{array}
$$
by which elements in $D_\triangle^\dagger$ are rational toric divisors for $\overline{D}^\dagger$. In particular, $\rho$ is expressed by $\rho = \sum_{\delta \in \Lambda^{(0)}} \frac{\rho^\delta}{m_\rho} e^{\overline{\rho} \dagger}$. Note that $\triangle^*_\rho$ can also be described by 
\be
\triangle^*_\rho = \{ y \in  L^*_{\RZ} \ | \ \langle y, \overline{\rho} \rangle \ \geq - \frac{\rho^\delta}{m_\rho} \  \ \mbox{for } \ \delta \in \Lambda^{(0)} \}. \tag{A10}
\ele(t*rho)
The homogeneous coordinates of $\PZ_{\Sigma}$ corresponding to the $\overline{D}$-structure are $\overline{z} = (\overline{z}_\delta)_{\delta} = \sum_{\delta \in \Lambda^{(0)} } \overline{z}_\delta e^{\overline{\delta} }$, related to $z = (z_\delta)_{\delta }$ in (\req(CZ)) by $\overline{z}_\delta = z^{m_\rho}_\delta$. Sections of an integral toric-divisor class of $\PZ_{\Sigma}$ are presented by $\overline{z}$-polynomials, and those for a rational toric-divisor class of $\PZ_{\Sigma}$ are expressed by polynomials of rational powers of $\overline{z}_\delta$'s. In particular, the rational toric-divisor classes induced by elements in $D^\dagger$ are presented by $z$-polynomials. Hence the formula (\req(sects)) still holds for sections of ${\cal O}(\iota^* \rho)$ over ${\PZ}_{\Sigma}$. The relation (\req(HX)) is valid for those quasi-smooth hypersurfaces $X$ defined by sections of integral toric-divisor classes. We summarize the obtained results as follows:
\begin{lem}\label{lem:sectK}
Let $\PZ_{\Sigma} (= \PZ_{\Sigma (\Lambda)}) $ be the ${\bf T} (L)$-toric variety defined by an integrable polytope $(\triangle, L)$ and a simplicial decomposition $\Lambda$ of $\partial \triangle$ with $\Lambda^{(0)} \subset L \cap \partial \triangle$. For $\delta \in \Lambda^{(0)}$, we denote by $m_\delta $ the positive integer characterized by  the primitive property of $\frac{\delta }{m_\delta }$ in $L$. Let ${\cal O}(\iota^* \rho)$ be the line bundle over ${\PZ}_{\Sigma}$ induced by an effective convex (rational) toric divisor $\rho = \sum_{\delta \in \Lambda^{(0)}} \rho^\delta e^{\delta \dagger} \in D_{\triangle}^\dagger$. Then 

{\rm (I)} The space of global sections of ${\cal O}(\iota^* \rho)$ over ${\PZ}_{\Sigma}$ has a combinatorial basis in (\req(sects)) with the coordinate representatives in (\req(zv*)).

{\rm (II)} If $\rho$ is an integrable toric divisor, i.e. $\frac{\rho^\delta }{m_\delta } \in {\rm \ZZ}$ for $\delta  \in \Lambda^{(0)}$, and $X$ is a quasi-smooth hypersurface of $\PZ_{\Sigma}$ defined by a section of ${\cal O} ( \iota^* \rho )$, then the space of sections of canonical bundle of $X$ has a basis in (\req(HX)).
\end{lem}  $\Box$ \par  \noindent
{\bf Remark.} (1) When $n=1$ in the above Lemma (II), $\PZ_{\Sigma} = \PZ^1$ with $L= \ZZ$ and $\frac{\delta }{m_\delta } = \pm 1$, in which case elements in (\req(sects))  are described by  monomials $z^k := z_1^k z_2^{d-k} ~ (0 < k < d)$, where $z= (z_1, z_2)$ is the homogeneous coordinates of $\PZ^1$, and $d=(\frac{\rho^0}{m_{\delta_0}}+\frac{\rho^1}{m_{\delta_1}})$. Then $X$ consists of $d$ general elements $x_j$ in $\PZ^1$, and $H^0(X)$ is a $d$-dimensional vector space with the basis $x_j^*$'s dual to $x_j$'s for $ 1 \leq j \leq d$. Then $H^{0, 0} (X)$ is the  $(d-1)$-dimensional subspace of $H^0(X)$ with $z^k$ corresponding to $\sum_{j=1}^d e^{2 \pi {\rm i}j k/d} x_j^*$ in (\req(HX)) for $(0 < k < d)$.
\par \noindent
(2) A section ${\sf w}^*$ in (\req(HX)) defines the effective divisor $\rho_{{\sf w}^*}= \sum_{\delta \in \Lambda^{(0)}} {\sf w}^*_\delta e^{\delta \dagger}$, which differs from $\rho$ by $\iota^*(\ell^*)$ for some $\ell^* \in L^*$, again satisfying the convex property.  The $(\rho_{{\sf w}^*})$-dual polytope $\triangle^*_{\rho_{{\sf w}^*}}$ is the translation of $\triangle^*_\rho$:  $\triangle^*_{\rho_{{\sf w}^*}} = \triangle^*_\rho-{\sf w}^*$.
\par \noindent
(3) By the discussion in \cite{R96},  the relations (\req(sects)) and (\req(HX)) are also valid for a convex toric divisor with $\rho^\delta \geq 0$. In this situation, the origin may not be in the interior of  the polytope $\triangle^*_\rho$. Furthermore, $\triangle^*_\rho$ could be a polytope of dimension less than $n$, in which case we set ${\rm Int} (\triangle^*_{\rho}) := \emptyset$ in (\req(HX)).
$\Box$ \par \vspace{.2in} \noindent
Next we extend (\req(sects)) and (\req(HX)) to formulas held for a face of $\triangle^*_\rho$.   For a $(n-m-1)$-face $F$ of $\triangle$ ($0 \leq m \leq n-1$), $L \cap F$ spans a $(n-m)$-dimensional subspace of $L_{\QZ}$, whose  intersection with $L$ gives rise to a $(n-m)$-sublattice $L_F^\perp$ of $L$. Let $L_F (:= L/L_F^\perp)$ be the quotient lattice of rank $m$ with the exact sequence of lattices:
\be
0 \longrightarrow L_F^\perp  \hookrightarrow L \stackrel{\wp}{\longrightarrow} L_F \longrightarrow 0 ,  \tag{A11}
\ele(LF)
where $\wp$ is the natural projection. A simplex ${\tt s} $ of a triangulation $\Lambda$ of $\partial \triangle$ whose interior is contained in the interior of $F$ must be in the $(n-m_{\tt s}-1)$-skeleton $\Lambda^{(n-m_{\tt s}-1)}$ with $m \leq m_{\tt s} \leq n-1$.  We define  
$$ 
{\tt s}^* =  \{ y \in  \triangle^*_\rho \ | \ \langle y, \delta \rangle  = -\rho^\delta  ~  {\rm ~ for} ~ \delta \in {\tt s} \cap \Lambda^{(0)} \},
$$
which is a $m$-face $F_\rho^*$ in $\partial \triangle^*_\rho$ with its dual $(n-m-1)$-$\rho$-face $F_\rho$ contained in $F$. The $(n-m_{\tt s})$-dimensional subspace of $L_{\QZ}$ spanned by one-parameter subgroups in the cone $\sigma ({\tt s})$, i.e. $L \cap \sigma ({\tt s})$, intersects $L$ on a $(n-m_{\tt s})$-sublattice $L_{\tt s}^\perp$,  whose quotient $L_{\tt s} (:= L/L_{\tt s}^\perp)$ in $L$ is a $m_{\tt s}$-lattice:
\bea(lll)
0 \longrightarrow L_{\tt s}^\perp  \hookrightarrow L \stackrel{\wp_{\tt s}}{\longrightarrow} L_{\tt s} \longrightarrow 0 , &  {\rm with} ~ L_{\tt s}^\perp  \subseteq L_F^\perp , & L_{\tt s}/L_F \simeq L_F^\perp/ L_{\tt s}^\perp .  \tag{A12}
\elea(Ls)
The ${\bf T} (L)$-orbit ${\rm O}_{\tt s}$ fixed by $ L \cap \sigma ({\tt s})$ is isomorphic to a $m_{\tt s}$-torus, whose closure $\overline{\rm O}_{\tt s}$  is a subvariety of $\PZ_{\Sigma}$. First we describe the toric structure of $\overline{\rm O}_{\tt s}$. Indeed, $\overline{\rm O}_{\tt s} = \cup_{{\tt s} \subset {\tt s}' \in \Lambda } {\rm O}_{s'} $. For a ${\tt s}' \in \Lambda^{(n-m_{\tt s} + k)}$ with ${\tt s} \subset {\tt s}' ~ (k \geq 0)$, the vertices of ${\tt s}'$ outside ${\tt s}$ form a $k$-simplex ${\tt s}'' \in \Lambda^{(k)}$. Denote 
$$ 
{\rm Star}({\tt s})= \bigcup_{k=0}^{m_{\tt s}-1} \{ {\tt s}'' \in \Lambda^{(k)} | {\tt s} \ast {\tt s}'' \in \Lambda^{(n-m_{\tt s} + k)} \}. 
$$
Under the projection $\wp_{\tt s}$ in (\req(Ls)), $\wp_{\tt s}({\tt s})$ is the origin $\overline{0}$ of $(L_{\tt s})_{\RZ}$, and $\wp_{\tt s} ({\rm Star}({\tt s}))$ becomes a  star of $\overline{0}$ in $(L_{\tt s})_{\RZ}$. Then  $\overline{{\tt s}''} ~ (:=\wp_{\tt s} ({\tt s}''))$  is a $k$-simplex with vertices in $L_{\tt s}$ for ${\tt s}'' \in {\rm Star}({\tt s})^{(k)}$. The convex hull of all $\overline{{\tt s}''}$, denoted by $\triangle ({\tt s})$, is a $m_{\tt s}$-polytope in $(L_{\tt s})_{\RZ}$ containing the origin so that $\wp_{\tt s} ({\rm Star}({\tt s}))$ form a triangulation $\Lambda ({\tt s})$ of the boundary  $\partial \triangle ({\tt s})$. Through $\wp_{\tt s}$ in (\req(LF)), the ${\bf T} (L)$-space $\overline{\rm O}_{\tt s}$ becomes a ${\bf T} (L_{\tt s})$-toric variety:
\be
\overline{\rm O}_{\tt s} = {\bf T} (L_{\tt s}){\rm -toric ~ variety ~ } \PZ_{\Sigma(\Lambda ({\tt s}))}, ~ ~ \Lambda ({\tt s}): {\rm triangulation ~ of} ~ \partial \triangle ({\tt s}).   \tag{A13}
\ele(Os)
Then ${\cal O}(\iota^* \rho)$ over $\PZ_{\Sigma}$ restricting on $\overline{\rm O}_{\tt s}$ gives rise to a line bundle over $\PZ_{\Sigma(\Lambda_{\tt s})}$, whose divisor class we are now going to derive. Let $F_\rho, F$ be the $(n-m-1)$-$\rho$-face and $(n-m-1)$-face of $\triangle$  which contain the interior of ${\tt s}$. Since ${\tt s}^* =  F^*_\rho $, using the coordinate form (\req(zv*)), one finds that the section corresponding to ${\sf v}^*$ in (\req(sects)) vanish on $\overline{\rm O}_{\tt s}$ when ${\sf v}^* \notin F^*_\rho $\footnote{Indeed, the conclusion can also be derived from the relation between the $\rho$-graph cone ${\cal C}_\rho$ and its dual cone ${\cal C}^*_\rho$ in the discussion of \cite{R96}  Section 3.}. So one needs only to consider those ${\cal O}(\iota^* \rho)$-sections corresponding to elements in $ F_\rho^* \cap L^*$. Note that vertices in $\Lambda ({\tt s})^{(0)}$ are of the form $\overline{\delta}''$ for $\delta'' \in {\rm Star}({\tt s})^{(0)}$, equivalently ${\tt s} \ast \delta''$ is a $(n-m)$-simplex in $\Lambda^{(n-m)}$. Hence $({\tt s} \ast \delta'')^*$ is a facet of $F_\rho^*$ with the expression $(s \ast \delta'')^* = \{ y \in F_\rho^* |  \langle y, \delta'' \rangle  = -\rho^{\delta''} \}$. The facets of $F_\rho^*$ are characterized  by $(m-1)$-faces $F_\rho^{\prime \prime *}$ with $F_\rho^{\prime \prime *} \subset F_\rho^{ *} $, equivalently, $(n-m)$-$\rho$-faces $F_\rho^{\prime \prime}$ with $F_\rho^{\prime \prime} \supset F_\rho $.  Since  $({\rm Star}({\tt s})^{(0)} \setminus F_\rho ) \cap F_\rho^{\prime \prime} \neq \emptyset$ when $F_\rho^{\prime \prime} \supset F_\rho $, we obtain
\be
F_\rho^*= \{ y \in L^*_{\RZ} | \langle  y, \delta \rangle = - \rho^\delta  ~ (\delta \in (F_\rho \cap \Lambda^{(0)}), ~ ~ \langle  y, \delta'' \rangle \geq - \rho^{\delta''} ~ (\delta'' \in {\rm Star}({\tt s})^{(0)}\setminus F_\rho) \}.  \tag{A14}
\ele(Freq)
If $ F_\rho^* \cap L^* = \emptyset$, the zeros of every section in (\req(sects)) contain $\overline{\rm O}_{\tt s}$. Otherwise, we
choose an element ${\rm v}_0^*$  in $ F_\rho^* \cap L^*$ as the base element. For convenience, ${\rm v}_0^*$ will be chosen in the interior of $F_\rho^*$ if $ {\rm Int} (F^*_{\rho}) \cap L^* \neq \emptyset$. 
Consider the $m$-dimensional subspace $\langle F_\rho^* - {\rm v}_0^* \rangle_{\QZ}$ of $L^*_{\QZ}$ spanned by $F_\rho^* - {\rm v}_0^*$. There is the $m$-lattice $L^*_{(F_\rho^*, {\rm v}_0^*)} = L^* \cap \langle F_\rho^* - {\rm v}_0^* \rangle_{\QZ}$ so that $ L^*_{(F_\rho^*, {\rm v}_0^*) ~ \QZ} = \langle F_\rho^* - {\rm v}_0^* \rangle_{\QZ}$. Since $L^*_{(F_\rho^*, {\rm v}_0^*)}$ vanishes on $F_\rho$, hence on $L_F^\perp$, one can identify $L^*_{(F_\rho^*, {\rm v}_0^*)}$ with the dual lattice $L_F^*$ of $L_F$ in (\req(LF)). Then $F_\rho^* - {\rm v}_0^*$ corresponds to an 
$m$-polytope in $L_{F ~ \RZ}^*$, denoted by $\triangle^\ddagger_{F_\rho}$:
\bea(lll)
\phi: & L_F^* ~ \stackrel{\sim}{\longrightarrow}   L^*_{(F_\rho^*, {\rm v}_0^*)} , & \triangle^\ddagger_{F_\rho} \longleftrightarrow F_\rho^* - {\rm v}_0^*  ,  \tag{A15}
\elea(tdag)
where  the isomorphism $\phi$ satisfies the relation $\langle  {\bf y} , \wp (x) \rangle = \langle  \phi ({\bf y} ), x \rangle$ for ${\bf y} \in L_F^*$ and $x \in L$. The integral elements in $\triangle^\ddagger_{F_\rho}$ are given by ${\bf y}_{{\rm v}^*} ~ ({\rm v}^* \in  F_\rho^* \cap L^*)$, where ${\bf y}_{{\rm v}^*}$'s are defined by the relation $\phi ({\bf y}_{{\rm v}^*}) = {\rm v}^* - {\rm v}_0^*$. Define the divisor $\overline{\rho}$ of the ${\bf T}(L_{\tt s})$-toric variety $\PZ_{\Sigma(\Lambda ({\tt s}))}$ in (\req(Os)): 
\bea(ll)
\overline{\rho} = \sum_{\delta'' \in {\rm Star}({\tt s})^{(0)}} \overline{\rho}^{\overline{\delta}''} e^{\overline{\delta}'' \dagger} \in D^\dagger_{\triangle ({\tt s})}, &
\overline{\rho}^{\overline{\delta}''}:= \rho^{\delta''}+ \langle {\rm v}_0^*, \delta''\rangle \in \ZZ_{\geq 0}~ ~{\rm for} ~ \delta'' \in {\rm Star}({\tt s})^{(0)}.   \tag{A16}
\elea(rhob)
Note that $\overline{\rho}^{\overline{\delta}''}=0$ when $\delta'' \in {\rm Star}({\tt s})^{(0)} \cap F_\rho$. 
By (\req(Freq)) and (\req(tdag)), one finds the equality of polytopes in $L_{F ~ \RZ}^*$:
\bea(ll)
\triangle^\ddagger_{F_\rho} = \triangle ({\tt s})^*_{\overline{\rho}}= \{ {\bf y} \in L^*_{F ~ \RZ} | \langle  {\bf y}, \overline{\delta}'' \rangle \geq - \overline{\rho}^{\delta''} ~ {\rm for} ~ \delta'' \in {\rm Star}({\tt s})^{(0)}) \}.  \tag{A17}
\elea(tdf)
The convex property of $F_\rho^*$ implies the convex condition of the divisor $\overline{\rho}$, which is effective when $ {\rm Int} (F^*_{\rho}) \cap L^* \neq \emptyset$. Using (\req(tdag)), (\req(rhob)), (\req(tdf)) and Lemma \ref{lem:sectK} (I), together with Remark (3) of Lemma \ref{lem:sectK}, one obtains the following results:
\begin{prop}\label{prop:srest}
Let ${\tt s} \in \Lambda^{(n-m_{\tt s}-1)}$ be a $(n-m_{\tt s}-1)$-simplex whose interior is contained in interior of a $(n-m-1)$-$\rho$-face $F_\rho$ and $(n-m-1)$-face $F$ (with $F_\rho \subseteq F$ and $0 \leq m \leq m_{\tt s} \leq n-1$). 
Denote by $\overline{\rm O}_{\tt s}$ the closure of the ${\bf T} (L)$-orbit associated to ${\tt s}$ in  $\PZ_{\Sigma (\Lambda)}$. 

{\rm (I)} $\overline{\rm O}_{\tt s}$ is isomorphic to the $m_{\tt s}$-dimensional ${\bf T} (L_{\tt s})$-toric variety ${\rm \PZ}_{\Sigma(\Lambda ({\tt s}))}$ in (\req(Os)). 

{\rm (II)} If the zero-loci of a generic section in (\req(sects)) does not contain $\overline{\rm O}_{\tt s}$, the restriction of ${\cal O} ( \iota^* \rho )$ on $\overline{\rm O}_{\tt s}$ is equivalent to the divisor class ${\cal O} ( \overline{\iota}^* \overline{\rho} )$ over ${\rm \PZ}_{\Sigma(\Lambda ({\tt s}))}$, where $\overline{\rho}$ is defined in (\req(rhob)) using a base element ${\rm v}_0^* \in  F_\rho^* \cap L^*$. The space of sections of ${\cal O} ( \overline{\iota}^* \overline{\rho} )$ over ${\rm \PZ}_{\Sigma(\Lambda ({\tt s}))}$ is described by 
\bea(lll) 
\Gamma ( {\rm \PZ}_{\Sigma(\Lambda ({\tt s}))} , {\cal O} ( \overline{\iota}^* \overline{\rho} )) &  \simeq & \bigoplus \{ {\rm  \CZ} {\sf v}^* \ | \ 
{\sf v}^* \in  F_\rho^* \cap L^* \ \}   \tag{A18}
\elea(Fsects)
so that the section with zero-divisor $ \overline{\rho}$ corresponds to ${\rm v}_0^*$.
\end{prop} $\Box$ \par  \noindent
Next we explore the relation between $L$-integral points in ${\rm Int}(F_\rho^*)$ and hypersurfaces of  $\overline{\rm O}_{\tt s}$ under the condition (II) of Proposition \ref{prop:srest}. Suppose  $\overline{\rho}$ in Proposition \ref{prop:srest} (II) is an integral toric divisor of ${\rm \PZ}_{\Sigma(\Lambda ({\tt s}))}$\footnote{In general, $\overline{\rho}$ is a rational, not integral, toric divisor in ${\rm \PZ}_{\Sigma(\Lambda ({\tt s}))}$ even when $\rho$ is an integral divisor of ${\rm \PZ}_{\Sigma}$.}, and $\overline{X}_{\tt s}$ is a quasi-smooth hypersurface in $\overline{\rm O}_{\tt s}$ defined by a generic section of ${\cal O} ( \overline{\iota}^* \overline{\rho} )$. First we consider the case when $m_{\tt s} =m$, i.e. ${\tt s} \in \Lambda^{(n-m-1)}$ with $ {\tt s} \subseteq F_\rho  \subseteq F$ for a $(n-m-1)$-$\rho$-face $F_\rho$ and $(n-m-1)$-face $F$. In this situation, by (\req(domi)) we have $L^\perp_{\tt s} = L^\perp_F$ and $L_{\tt s} = L_F$. By Lemma \ref{lem:sectK} (II), together with Remark (3) of Lemma \ref{lem:sectK}, one finds the combinatorial representation of canonical forms of $\overline{X}_{\tt s}$:
\bea(lll)
H^{m-1, 0} (\overline{X}_{\tt s}) & \simeq  & \bigoplus \{ {\rm \CZ} {\sf w}^* \ | \ {\sf w}^* \in  {\rm Int} (F_\rho^*) \cap L^* \ \}. \tag{A19}
\elea(KXb)
Note that in the above formula, the data on the right are the same for all $\overline{X}_{\tt s}$ whenever ${\tt s} \in \Lambda^{(n-m-1)}$ with $s \subseteq F_\rho$.  Indeed under the dominating morphism (\req(domi)), $\varphi (\overline{\rm O}_{\tt s})$ is equal to the closure $\overline{\bf O}_{F_\rho}$ of ${\bf T}(L)$-orbit ${\bf O}_{F_\rho}$ in the $\rho$-minimal toric variety $\PZ_{\Sigma_{\rho ~ 0}}$. With the same argument as the toric structure  of $\overline{\rm O}_{\tt s}$ in (\req(Os)), one finds $
\overline{\bf O}_{F_\rho}$ is isomorphic to a complete ${\bf T} (L_F)$-toric variety $\PZ_{\Sigma(\Lambda ({F_\rho}))}$,
where $L_F$ is the lattice in (\req(LF)), and $\Lambda ({F_\rho})$ is a polytope decomposition of the boundary of a $m$-polytope in $L_{F~\RZ}$ obtained by the $\wp$-projection of all $\rho$-faces $F'_\rho \supset F_\rho$. Furthermore, as the relation between $\rho$ and $\rho_0$ in (\req(domi)), the divisor $\overline{\rho}$ in $\overline{\rm O}_{\tt s}$ is the $\varphi$-pull-back of a toric divisor $\overline{\rho}_0$ in $\overline{\bf O}_{F_\rho}$: ${\cal O} (\overline{\iota}^* \overline{\rho}) = \varphi_{\tt s}^* {\cal O} (\overline{\iota}_0^* \overline{\rho}_0)$. Here $\varphi_{\tt s}$ denotes the restriction of $\varphi$ on $\overline{\rm O}_{\tt s}$. Hence $\Gamma ( {\rm \PZ}_{\Sigma(\Lambda (F_\rho))} , {\cal O} ( \overline{\iota}_0^* \overline{\rho}_0 ))$ is isomorphic to vector spaces in (\req(Fsects)), and $\overline{X}_{\tt s}$ in $\overline{\rm O}_{\tt s}$ is indeed the pull-back of a hypersurface $\overline{X}_{F_\rho}$ in $\overline{\bf O}_{F_\rho}$: $\overline{X}_{\tt s} = \varphi^*(\overline{X}_{F_\rho})$. Therefore the combinatorial data in (\req(KXb)) can be regarded as those for $\overline{X}_{F_\rho}$.  A similar consideration can be carried over to the general case for $m_{\tt s} \geq m$ in Proposition \ref{prop:srest} (II). Indeed,  $\varphi (\overline{\rm O}_{\tt s}) = \overline{\bf O}_{F_\rho}$, equivalently the restriction of (\req(domi)) defines a surjective morphism 
\bea(ll)
\varphi_{\tt s}: \PZ_{\Sigma(\Lambda ({\tt s}))} \longrightarrow \PZ_{\Sigma(\Lambda ({F_\rho}))} , & {\cal O} (\overline{\iota}^* \overline{\rho}) = \varphi_{\tt s}^* {\cal O} (\overline{\iota}_0^* \overline{\rho}_0)  \tag{A20}
\elea(vphis)
so that $\overline{X}_{\tt s}$  is the pull-back of a hypersurface $\overline{X}_{F_\rho}$ in $\PZ_{\Sigma(\Lambda (F_\rho)}$ defined by zeros of a  ${\cal O} (\overline{\iota}_0^* \overline{\rho}_0)$-section: $\overline{X}_{\tt s} = \varphi_{\tt s}^{-1}(\overline{X}_{F_\rho})$.
By the toric data associated to the fibration (\req(vphis)), the generic fiber ${\sf f}_{\tt s}$ of $\varphi_{\tt s}$ is isomorphic a complete ${\bf T} (L^\perp_F/L^\perp_{\tt s})$-toric variety $\PZ_{\Sigma(\Lambda ({\tt s}, F_\rho))}$, where $L^\perp_F, L^\perp_{\tt s}$ are lattices in (\req(Ls)), and $\Lambda ({\tt s}, F_\rho))$ is a triangulation of the boundary of a $(m_{\tt s}-m)$-polytope in 
$L^\perp_{F ~ \RZ}/L^\perp_{{\tt s}~\RZ}$ generated by ${\rm Star}({\tt s})^{(0)} \cap F_\rho$. Hence $\overline{X}_{\tt s}$ is a fibration over  $\overline{X}_{F_\rho}$ with the general fiber ${\sf f}_{\tt s}$. Since $H^{j, 0}({\sf f}_{\tt s}) = 0 $ for $1 \leq  j \leq m_{\tt s}-m$, $H^{k, 0}(\overline{X}_{\tt s})$ vanishes for $m \leq k \leq m_{\tt s}-1$ and $(m-1, 0)$-forms  of $\overline{X}_{\tt s}$ are $\varphi_{\tt s}$-pull-back of canonical forms of $\overline{X}_{F_\rho}$  given by (\req(KXb)). We summarize the result as follows: 
\begin{prop}\label{prop:sform}
Let a simplex ${\tt s} \in \Lambda^{(n-m_{\tt s}-1)}$, $(n-m-1)$-$\rho$-face $F_\rho$, $(n-m-1)$-face $F$, $\overline{\rm O}_{\tt s} = {\rm \PZ}_{\Sigma(\Lambda ({\tt s}))}$, and $\overline{\rho}$ be the same as in Proposition \ref{prop:srest} {\rm (II)}. Suppose  $\overline{\rho}$ is an integral toric divisor of ${\rm \PZ}_{\Sigma(\Lambda ({\tt s}))}$, and $\overline{X}_{\tt s}$ is a quasi-smooth hypersurface in $\overline{\rm O}_{\tt s}$ defined by a generic section of ${\cal O} ( \overline{\iota}^* \overline{\rho} )$. Then $H^{k, 0}(\overline{X}_{\tt s}) =0$ for $k \geq m$, and $H^{m-1, 0}(\overline{X}_{\tt s})$ is isomorphic to the space in (\req(KXb)). 
\end{prop}
$\Box$ \par  \noindent
{\bf Remark.} Indeed, the $\overline{X}_{\tt s}$ in the above proposition is a fibration over a quasi-smooth hypersurface $\overline{X}_{F_\rho}$ of $\PZ_{\Sigma(\Lambda (F_\rho)}$ in (\req(vphis)) with the general fiber ${\sf f}_{\tt s}$ being a complete toric variety of dimension $m_{\tt s}-m$. The elements in (\req(KXb)) are induced from the canonical forms of the base $\overline{X}_{F_\rho}$.
$\Box$ \par \vspace{.1in} \noindent
We now discuss the relation between the homogeneous coordinates of ${\PZ}_{\Sigma}$ in (\req(prin)) and $\overline{\rm O}_{\tt s}$ in (\req(Os)):
\be
\pi_{\tt s} : \CZ_{\widetilde{\Sigma}(\Lambda ({\tt s}))} \longrightarrow {\PZ}_{{\Sigma(\Lambda ({\tt s}))}}  ~ ~ ( = \overline{\rm O}_{\tt s}) .  \tag{A21}
\ele(prins)
As in (\req(exa)) and (\req(exa*)), we have the exact sequences
\bea(ll)
0 \longrightarrow {\bf n}_{\triangle ({\tt s})} \stackrel{\iota_{\tt s}}{\longrightarrow} D_{\triangle ({\tt s})} \stackrel{\beta_{\tt s}}{\longrightarrow} L_{F ~ 0} \longrightarrow  0 , & 
0 \longrightarrow   L_{F ~ 0}^* \stackrel{ \beta_{\tt s}^*}{\longrightarrow} D_{\triangle ({\tt s})}^*
\stackrel{\iota_{\tt s}^*  }{\longrightarrow}  {\bf n}_{\triangle ({\tt s})}^*  \longrightarrow 0  ,  \tag{A22}
\elea(exas)
where $L_{F ~ 0}$ is the $m$-sublattice of $L_F$ generated by $\Lambda ({\tt s})^{(0)}$, and $\beta_{\tt s} ( e^{\overline{\delta}''} ) := \overline{\delta}''$ for $\overline{\delta}'' \in \Lambda ({\tt s})^{(0)}$. By (\req(zv*)), (\req(tdag)) and (\req(rhob)),  ${\sf v}^*$ in (\req(Fsects)) can be expressed in terms of the homogeneous coordinates $\overline{z} = (\overline{z}_{\overline{\delta}''})_{\delta'' \in {\rm Star}({\tt s})^{(0)}}$ of (\req(prins)) as follows:
\be
{\overline{z}}^{{\sf v}^*} = \prod_{\delta'' \in {\rm Star}({\tt s})^{(0)}} \overline{z}_{\overline{\delta}''}^{{\sf v}_{\overline{\delta}''}^*},  ~ ~ ~ {\sf v}^*_{\overline{\delta}''} (:= \langle {\bf y}_{{\rm v}^*}, \overline{\delta}'' \rangle_F + \rho^{\overline{\delta}''}) = {\sf v}^*_{\delta''} ~ ~ {\rm for} ~ \delta'' \in {\rm Star}({\tt s})^{(0)}. \tag{A23}
\ele(zv*s)
On the other hand, the coordinates (\req(zv*)) of ${\sf v}^* \in  F_\rho^* \cap L^*$ in (\req(sects)) are given by 
\be
z^{{\sf v}^*} = \big(\prod_{\delta'' \in {\rm Star}({\tt s})^{(0)}} z_{\delta'' }^{{\sf v}^*_{\delta''} } \big) 
\big( \prod_{\delta \notin (F_\rho \cup {\rm Star}({\tt s})^{(0)})} z_{\delta }^{{\sf v}^*_\delta } \big)  ~ ~ {\rm for} ~ {\sf v}^* \in  F_\rho^* \cap L^* \tag{A24}
\ele(zFrho)
since ${\sf v}^*_\delta = 0$ for $\delta \in \Lambda^{(0)} \cap F_\rho$. The coordinates (\req(zv*s)) are obtained from (\req(zFrho)) by setting $z_\delta=1$ for $\delta \notin {\rm Star}({\tt s})^{(0)}$, which gives rise to the regular embedding from (\req(prins)) into (\req(prin)):
$$
\begin{array}{lll}
\CZ_{\widetilde{\Sigma}(\Lambda ({\tt s}))} \longrightarrow \CZ_{\widetilde{\Sigma}}, & \overline{z} = (\overline{z}_{\overline{\delta}''})_{\delta'' \in {\rm Star}({\tt s})^{(0)}} \mapsto
z = (z_\delta)_{\delta \in \Lambda^{(0)} } & z_\delta = \left\{ \begin{array}{ll}\overline{z}_{\overline{\delta}} & {\rm if} ~ \delta \in {\rm Star}({\tt s})^{(0)} \\
1 & {\rm otherwise} ,
\end{array} \right.
\end{array}
$$
induced by the divisor-embedding of (\req(exas)) into  (\req(exa)):
$$
 D_{\triangle ({\tt s})} \hookrightarrow D_{\triangle} , ~ ~ \overline{\delta}'' \mapsto {\delta}'' ~ ~ {\rm for} ~ \delta'' \in {\rm Star}({\tt s})^{(0)}.
$$
Hence the lattice ${\bf n}_{\triangle ({\tt s})}$ in (\req(exas)) is given by ${\bf n}_{\triangle ({\tt s})} = {\rm proj}_{\triangle ({\tt s})}(\beta^{-1}(L_F^\perp))$, where $\beta, L_F^\perp$ are defined in (\req(exa)), (\req(LF)) respectively, and ${\rm proj}_{\triangle ({\tt s})}$ is the projection of $D_{\triangle}$ onto $D_{\triangle ({\tt s})}$.

Consider the case when the $\rho$-dual polytope $\triangle^*_\rho$ is $L^*$-integral. Note that in the case $\rho = - \kappa$, the integral $\triangle^*_{-\kappa}$ is the same as the reflexive polytope condition of $(\triangle, L)$.  By (\req(t*rho)), conditions in Lemma \ref{lem:sectK} (II), Proposition \ref{prop:srest} (II), and Proposition \ref{prop:sform}  are all satisfied. Furthermore, by \cite{R96} Theorem 1 and 2\footnote{The convention $\Lambda$ in this paper is different from that in \cite{R96} Theorem 2, where $\Lambda$ denotes the polyhedral cone of the boundary of $\rho$-graph cone over $\Sigma$.}, a combinatorial basis for the Picard group of a generic hypersurface of some special toric variety $\PZ_{\Sigma}$ is described as follows:
\begin{prop}\label{prop:rhoRf}
Let $\Lambda$ be a simplicial decomposition of $\partial \triangle$ of an integral polytope $(\triangle, L)$ with $\Lambda^{(0)} \subset L \cap \partial \triangle$, and $\rho$ be an effective convex toric divisor in (\req(rho)). Suppose the $\rho$-dual polytope $\triangle^*_\rho$ is $L^*$-integral, and  $\Lambda$ satisfies the following condition
$$
\Lambda^{(0)} = L \ \bigcap \ \big( \partial \triangle - 
\bigcup \{ {\rm Int}(F_\rho)  \ | \  F_\rho :  {\rm codim-1} ~ \rho{\rm -face \
of \ } \triangle \} \big).
$$
Let ${\tt s} \in \Lambda^{(n-m_{\tt s}-1)}$ be a $(n-m_{\tt s}-1)$-simplex whose interior is contained in interior of a $(n-m-1)$-$\rho$-face $F_\rho$ and $(n-m-1)$-face $F$ (with $F_\rho \subseteq F$ and $0 \leq m \leq m_{\tt s} \leq n-1$), and $\overline{\rm O}_{\tt s}$ be the closure $\overline{\rm O}_{\tt s}$ of ${\bf T} (L)$-orbit associated to ${\tt s}$. Consider a hypersurface $X$ of ${\rm \PZ}_{\Sigma}$ defined by a generic section of ${\cal O}(\iota^* \rho)$. Then

{\rm (I)} $X$ is quasi-smooth,  and we have the following combinatorial representation for the Picard group of $X$: 
$$
{\rm Pic} ( X  )_{\rm \QZ} \simeq {\bf n}^\dagger_{\triangle~ {\rm \QZ}} 
\oplus \bigoplus_{F_\rho , \nu_{F_\rho} } {\rm \QZ} \nu_{F_\rho},
$$
where the index $F_\rho$ runs over the codimensional 2 $\rho$-faces of $\partial \triangle$, and $\nu_{F_\rho} \in ( {\rm Int} ( F_\rho ) \cap L ) \times ( {\rm Int} ( F_\rho^* ) \cap L^* )$.

{\rm (II)} $\overline{\rm O}_{\tt s}$ is isomorphic to a $m_{\tt s}$-dimensional toric variety  ${\rm \PZ}_{\Sigma (\Lambda)}$ in (\req(Os)). The restriction of ${\cal O} ( \iota^* \rho )$ on $\overline{\rm O}_{\tt s}$ is equivalent to the divisor class ${\cal O} ( \overline{\iota}^* \overline{\rho} )$ over ${\rm \PZ}_{\Sigma(\Lambda ({\tt s}))}$ where $\overline{\rho}$ is an integral toric divisor of ${\rm \PZ}_{\Sigma(\Lambda ({\tt s}))}$ defined in (\req(rhob)), and 
the space of sections of ${\cal O} ( \overline{\iota}^* \overline{\rho} )$ over ${\rm \PZ}_{\Sigma(\Lambda ({\tt s}))}$ has a basis represented by $L^*$-integral points of $F_\rho^*$ in (\req(Fsects)).

{\rm (III)} $X$ intersects $\overline{\rm O}_{\tt s}$ on a quasi-smooth hypersurface $\overline{X}_{\tt s}$ of $\PZ_{\Sigma (\Lambda)}$ with $H^{k, 0}(\overline{X}_{\tt s}) =0$ for $k \geq m $, and $H^{m-1, 0}(\overline{X}_{\tt s})$ isomorphic to the space in (\req(KXb)). 
\end{prop}
$\Box$ \par 
\end{document}